\documentclass[12pt]{amsart}
\usepackage{amscd,amssymb}
\theoremstyle{plain}
\newtheorem{prop}[subsection]{Proposition}
\newtheorem{thm}[subsection]{Theorem}
\newtheorem{lem}[subsection]{Lemma}
\newtheorem{cor}[subsection]{Corollary}

\theoremstyle{remark}
\newtheorem{rem}[subsection]{Remark}

\newtheorem*{conv}{Conventions}
\newtheorem*{claim}{Claim}

\theoremstyle{definition}
\newtheorem{exm}[subsection]{Example}

\numberwithin{equation}{section}

\newcommand{\A}{{\mathcal A}}
\newcommand{\calS}{{\mathcal S}}
\newcommand{\W}{{\mathcal W}}
\newcommand{\Z}{{\mathbb Z}}
\newcommand{\R}{{\mathbb R}}
\newcommand{\C}{{\mathbb C}}
\newcommand{\M}{{M'}}
\newcommand{\X}{{X'}}
\newcommand{\D}{{\Delta }}
\newcommand{\p}{\partial }
\newcommand{\wab}{\widetilde{\ab}}
\newcommand{\FTC}{fundamental formula of Fox Calculus}
\newcommand{\DS}{\displaystyle }
\newcommand{\hB}{\widehat B}
\newcommand{\J}{\widehat J}
\newcommand{\hL}{\widehat {\Lambda }}
\newcommand{\hD}{\widehat {\Delta }}
\newcommand{\LT}{\text {\scshape lt}}
\newcommand{\G}{{\mathcal G}}
\newcommand{\locB}{{B^{\operatorname{cc}}}}
\newcommand{\loc}{{\operatorname{cc}}}
\newcommand{\hlocB}{{\widehat{B}^{\operatorname{cc}}}}

\renewcommand{\a}{{\alpha }}
\renewcommand{\b}{{\beta }}
\renewcommand{\c}{{\gamma }}
\renewcommand{\d}{{\delta }}
\newcommand{\s}{{\sigma }}
\renewcommand{\L}{{\Lambda }}

\renewcommand{\ll}{{\ell }}
\renewcommand{\H}{{\mathcal H}}

\DeclareMathOperator{\rank}{rank}
\DeclareMathOperator{\im}{im}
\DeclareMathOperator{\coker}{coker}
\DeclareMathOperator{\id}{id}
\DeclareMathOperator{\Aut}{Aut}

\DeclareMathOperator{\IA}{IA}
\DeclareMathOperator{\CA}{CA}
\DeclareMathOperator{\GL}{GL}
\DeclareMathOperator{\ab}{ab}
\DeclareMathOperator{\Span}{span}
\DeclareMathOperator{\pr}{pr}
\DeclareMathOperator{\gr}{gr}

\begin{document}

\title[Alexander Invariants of Arrangements]%
{Alexander Invariants of Complex Hyperplane Arrangements}
\author{Daniel C.~Cohen}
\address{Department of Mathematics,
Louisiana State University,
Baton Rouge, LA 70803}
\email{cohen@math.lsu.edu}

\author{Alexander I.~Suciu}
\address{Department of Mathematics,
Northeastern University,
Boston, MA 02115}
\email{alexsuciu@neu.edu}

\thanks{The first author was partially supported by 
grant LEQSF(1996-99)-RD-A-04 from the Louisiana Board of Regents
and by a grant from the Louisiana State University Council on Research.}
\thanks{The second author was partially supported by N.S.F.~grant
DMS--9504833, and an RSDF~grant from Northeastern University.}



\begin{abstract}Let $\A $ be an arrangement of $n$ complex hyperplanes.
The fundamental group of the complement of $\A$ is determined by a braid
monodromy homomorphism, $\a :F_{s}\to P_{n}$. Using the Gassner
representation of the pure braid group, we find an explicit presentation 
for the Alexander invariant of $\A$.  From this presentation, we obtain 
combinatorial lower bounds for the ranks of the Chen groups of $\A$. 
We also provide a combinatorial criterion for when these lower 
bounds are attained.
\end{abstract}

\maketitle

\section*{Introduction}\label{sec:intro}

Let $\A =\{H_{1},\dots ,H_{n}\}$ be an arrangement of hyperplanes in 
$\C ^{d}$, with complement $M=\C ^{d} \setminus \cup_{i=1}^n H_{i}$, 
and group $G=\pi _{1}(M)$.  Let $\M$ be the maximal abelian cover, 
corresponding to the abelianization $\ab: G \to \Z^n$.  
The action of $\Z^n$ on $\M$ puts on $H_*(\M)$ the structure 
of a module over the group ring $\Z\Z^n$.  
This ring can be identified with the ring of Laurent polynomials 
$\L =\Z [t_{1}^{\pm 1},\dots ,t_{n}^{\pm 1}]$, with $t_i$ corresponding 
to a standardly oriented meridional loop around $H_i$.    
The object of our study is the Alexander invariant, $B(\A)=H_{1}(\M)$, 
viewed as a module over the ring $\L $.

Let $L(\A)$ denote the intersection lattice of $\A$, with rank function
given by codimension (see \cite{OT} as a general reference for arrangements).
Let $s$ denote the cardinality of $L_2(\A)$, the set of rank two elements
in $L(\A)$.  From the defining polynomial of $\A$, one can compute 
the Moishezon-Libgober braid monodromy homomorphism, $\a:F_s\to P_n$, 
see~\cite{CS3}.  This homomorphism determines a finite presentation for 
the group of the arrangement:  
$G=\langle t_{1},\dots ,t_{n}\mid \a _{k}(t_{i})=t_{i} \rangle$, 
where $\a_1,\dots ,\a_s$ generate the image of $\a$.  
The braid monodromy may also be used to 
obtain a finite presentation for the Alexander invariant $B(\A)$.  
We accomplish this here, by means of the Gassner representation, 
$\Theta :P_{n} \to \GL (n,\Lambda )$, the Fox calculus, 
and homological algebra.  

Surprisingly, the size of the presentation depends only on the first 
two betti numbers of the complement:  there are $\binom{n}{2}$ generators  
and $\binom{n}{3}+b_{2}(M)$ relations.  When $\A $ is the
complexification of a real arrangement, the presentation of $B(\A)$ can
be simplified to $\binom{n}{2} - b_{2}(M)$ generators and
$\binom{n}{3}$ relations.  More generally, if $G$ is the group of a
collection of $s$ basis-conjugating automorphisms of a finitely
generated free group $F_n$, our methods yield a presentation of $B(G)$
with $\binom{n}{2}$ generators and $\binom{n}{3}+ns$ relations.  In
particular, the Alexander invariant of any pure link has such a
presentation (with $s=1$).  This should be compared with the general
situation for links in $S^3$, where there is no  upper bound on the
number of relations, see~\cite{Ma}.   

Note that the Alexander invariant is isomorphic to $G'/G''$ 
(with the usual $G/G'$ action), and so depends only on the 
isomorphism type of $G$.  Consequently, we may obtain invariants of 
an arrangement $\A$ from the module $B=B(\A)$ and its presentation.
For instance, if $\Delta$ is a presentation matrix for $B$, 
the elementary ideal $E_k(B)$ is defined to be the ideal generated
by the codimension $k$ minors of $\Delta$.  It is well-known that these
ideals depend only on the module $B$.  These ideals, and the closely related
characteristic varieties, arise in the study of plane algebraic curves; 
see for instance the recent works of Hironaka~\cite{eko} and Libgober~\cite{L2}.
The structure of the elementary ideals and characteristic varieties of 
the Alexander invariant an arrangement will be the subject of a future work.

In this paper, we focus on Chen groups.  The Chen groups of $G$ 
are the lower central series quotients of the maximal metabelian 
quotient $G/G''$.  
Using an observation of Massey~\cite{Ma} relating the Chen groups and 
the Alexander invariant, together with Mora's tangent cone algorithm, 
we obtain an algorithm for computing the Chen groups of an 
arrangement from the presentation of the Alexander invariant $B$.
The ranks of the Chen groups often serve to distinguish the 
groups of combinatorially ``similar'' arrangements.  This is particularly
useful for fiber-type arrangements, where the ranks of the lower central
series quotients of $G$ itself are determined by the exponents
of the arrangement.  On the other hand, we know of no combinatorially
equivalent arrangements whose Chen groups differ.

The precise relation between the Chen groups and 
the intersection lattice of a central arrangement $\A$ is not known.  
We obtain partial results toward this end here.  
To each element $V\in L_2(\A)$, we associate
a ``local'' Alexander invariant $B_V$.  
Algebraic considerations yield a surjective homomorphism $B \to \locB$,
where $\locB=\oplus_V B_V$ is the ``coarse combinatorial Alexander
invariant'' of $\A$, determined by (only) the 
multiplicities of the elements of $L_2(\A)$.  
From this map, we obtain combinatorial lower bounds
on the ranks of the Chen groups of $\A$.  These ranks are determined by
the $I$-adic completion, $\hB$, of the Alexander invariant $B$, where $I$
is the augmentation ideal of $\Lambda$.  We find a 
combinatorial criterion for when the completion of the Alexander
invariant of $\A$ decomposes as a direct sum, i.e.,~
$\hB\xrightarrow{\sim}\hlocB$.   We also obtain a combinatorial 
formula for the rank of the third Chen group of any arrangement.

The above results may be viewed as evidence that the ranks, $\theta_k$, 
of the Chen groups of $\A$ are combinatorially determined.  
In~\cite{CS1}, we conjectured an explicit combinatorial formula
for $\theta_k$, for sufficiently large $k$.  
This formula involved the number $\beta$ of 
subarrangements of $\A$ 
lattice-isomorphic to the braid arrangement $\A_4~\subset~\C^4$.  
In the present context, we show by example that if $\beta>1$,
then the map $\hB \to \hlocB$ is not an isomorphism.  Other examples exhibit
combinatorially different ways this map can fail to be an isomorphism.
These provide counterexamples to the aforementioned formula, and illustrate
the subtlety of the relationship between the Chen groups and the lattice of
an arrangement.

Our results on Chen groups parallel a portion of Falk's work 
on the LCS quotients of an arrangement group.  The combinatorial 
lower bounds we obtain for the ranks of the Chen groups are 
analogous to those for the ranks, $\phi_k$, of the LCS quotients 
found in~\cite{F2}. Moreover, the formula we obtain 
for $\theta_3 = \phi_3$ may be viewed as dual to the description 
of $\phi_3$ found in~\cite{F1}, \cite{F2}.  
The precise relationship between the Chen groups and LCS quotients 
of an arrangement will be explored elsewhere.  

The structure of the paper is as follows:

\begin{itemize}
\item
In section~\ref{sec:AlexanderChen}, we review Alexander invariants 
and Chen groups, and present a Groebner basis algorithm for 
determining the latter.  The section concludes with an analysis 
of the Alexander invariant and Chen groups of a product of spaces.  
\item
In section~\ref{sec:FoxCalc}, we introduce our basic computational 
tools: the Fox free differential calculus and the Magnus representations.  
\item
In section~\ref{sec:AlexFreeAuto}, we study the Alexander invariant 
of the group of a free automorphism.  An explicit 
presentation is given when the automorphism is basis-conjugating. 
\item
In section~\ref{sec:LocalAlexInv}, we find presentations for the 
local Alexander invariants of an arrangement.  
\item
In section~\ref{sec:AlexInvArr}, the presentation for the Alexander 
invariant of an arrangement is obtained. 
\item
In section~\ref{sec:DecompAlex}, the homomorphism $B \to \locB$ is 
defined, and its completion proven to be an isomorphism when a 
certain criterion is satisfied.  
\item
In section~\ref{sec:Combinatorics}, the aforementioned criterion is
shown to be combinatorial, and  
lower bounds on the ranks
of the Chen groups of an arrangement are obtained.  
\item
In section~\ref{sec:Examples}, we illustrate 
our results
by means of several explicit examples.  
\end{itemize}

\begin{conv} Given a group $G$, we will denote by $\Aut(G)$ the 
group of {\em right} automorphisms of $G$, with multiplication 
$\alpha\cdot \beta=\beta\circ\alpha$.  We will regard all
modules over the group ring $\Z G$ as {\em left} modules.
Elements of the free module $(\Z G)^{n}$ are viewed as {\em row}
vectors, and $\Z G$-linear maps $(\Z G)^{n}\to (\Z G)^{m}$ are
viewed as $n\times m$ matrices which act on the {\em right}
(so that the matrix of $B\circ A$ is $A\cdot B$).  We will write 
$A^{\top }$ for the transpose of $A$, and 
$\left( A_1\ \cdots\ A_s \right)^{\top}$ for 
$\bigl(\begin{smallmatrix}A_1 \\ \dots \\
A_s \end{smallmatrix}\bigr)$.  
If $\phi :G\to H$ is a homomorphism, $\tilde {\phi }:\Z G\to \Z H$ 
denotes its $\Z$-linear extension to group rings.  
We will abuse notation and also write
$\tilde {\phi }:(\Z G)^{n}\to (\Z H)^{n}$ for the map
$\oplus _{1}^{n} \tilde \phi $.
\end{conv}

\section{Alexander Invariants, Chen Groups, and Products}
\label{sec:AlexanderChen}

We start by reviewing the definition of the Alexander 
invariant of a finite complex.  We then present an algorithm for computing 
the ranks of the Chen groups of a group, based on a presentation of 
this module.  Finally, we determine the structure of the 
Alexander invariant of a product of spaces in terms of 
those of the factors. 

\subsection{Alexander Invariants}\label{subsec:Alexander} 
Let $M$ be a path-connected space that has the homotopy type of a finite
CW-complex.  Let $G=\pi _{1}(M,*)$ be the fundamental group, and 
$K=H_{1}(M)$ its abelianization.  Let $p:\M \to M$ be the maximal abelian
cover. The action of $K$ on $\M $ passes to an action of $K$ on the homology
groups $H_*(\M)$. This defines on $H_*(\M)$ the structure of a module over
the group ring $\Z K$.  The $\Z K$-module $B=H_{1}(\M)$ is called the (first)
{\em Alexander invariant} of $M$.  Closely related to it is the (first) 
{\em Alexander module}, $A=H_{1}(\M , p^{-1}(*))$.  These two modules,
together with the augmentation ideal $I=IK=\ker (\epsilon :\Z K \to \Z )$,
comprise the Crowell exact sequence, $0\to B\to A \to I\to 0$, of \cite{Cr1}.

The two Alexander modules depend only on the group $G$.  Indeed, 
$A=\Z K \otimes_{\Z G} IG$, with $K=G/G'$ acting by multiplication on the 
left factor, and $B=G'/G''$, with the action of $K$ defined by 
the extension $1\to G'/G''\to G/G'' \to G/G'\to 1$.  Since $M$ is by 
assumption a finite complex, $G$ is a finitely presented group.  
Hence, the $\Z K$-module $A$ is finitely presented; 
Fox's free differential calculus 
provides an explicit presentation
(see \cite{Cr2}, and also sections \ref{sec:FoxCalc} and 
\ref{sec:AlexFreeAuto}).  Less evident, but still true, is 
the fact that $B$ also admits a finite presentation as a $\Z K$-module
(see \cite{Cr1}, \cite{Ma}, and also section \ref{sec:AlexFreeAuto}).  

\subsection{Chen Groups}\label{subsec:Chen}
Let $\Gamma _{k}(G)$ denote the $k^{\text{th}}$ lower central series 
subgroup of $G$, defined inductively by $\Gamma _{1}(G)=G$ and 
$\Gamma _{k+1}(G)=[\Gamma _{k}(G), G]$ for $k \ge 1$.  
The projection of $G$ onto its maximal metabelian quotient $G/G''$ 
induces an epimorphism
\begin{equation*} \label{eq:Gamma}
{\frac{\Gamma _{k}(G)}{\Gamma _{k+1}(G)}}
\twoheadrightarrow{\frac{\Gamma _{k}(G/G'')}
{\Gamma _{k+1}(G/G'')}}
\end{equation*}
from the $k^{\text{th}}$ lower central series quotient of $G$ to the
$k^{\text{th}}$ Chen group of $G$.  Since $G$ is finitely presented, 
these quotients are finitely generated abelian groups, whose ranks we 
will denote by $\phi_k$, respectively $\theta_k$.  It is readily seen 
that $\phi_k=\theta_k$ for $k\le 3$, and $\phi_k \ge \theta_k$ for $k> 3$.

The Chen groups of $G$ can be determined from the Alexander invariant of $G$.  
Indeed, Massey \cite{Ma} noted the following isomorphism, for $k\ge 2$:
\begin{equation*} \label{eq:massey}
{\frac{\Gamma _{k}(G/G'')}{\Gamma _{k+1}(G/G'')}}=
{\frac{I^{k-2} B}{I^{k-1}B}}.
\end{equation*}
Thus, the Chen groups are determined by 
$\gr B=\bigoplus _{k\ge 0} I^{k}B/I^{k+1}B$, viewed as a graded module 
over the graded ring $\gr \Z K =\bigoplus _{k\ge 0} I^{k}/I^{k+1}$.  

Now assume $K$ is free abelian, and fix a system of generators, 
$t_{1},\dots , t_{n}$.   The group ring $\Z K$ can be identified 
with the ring of Laurent polynomials in $n$ variables, 
$\L =\Z [t_{1}^{\pm 1},\dots ,t_{n}^{\pm 1}]$.  The ring $\L$ 
can be viewed as a subring of the formal power series ring 
$P=\Z [[x_{1},\dots ,x_{n}]]$ via the ``Magnus embedding,'' given by
$t_{i}\mapsto 1-x_{i}$ and $t_{i}^{-1}\mapsto \sum _{k=0}^{\infty }x_{i}^{k}$.
Let $\hL =\varprojlim \L /I^{k}$ be the completion of $\L $ relative to the 
$I$-adic topology.  Then, the Magnus embedding extends to a ring isomorphism
$\hL \xrightarrow{\sim} P$.

Consider the $\mathfrak{m}$-adic filtration on $P$, where 
$\mathfrak{m}=\langle x_{1},\dots , x_{n} \rangle $, and its associated graded 
ring, $\gr P=\bigoplus _{k\ge 0}\mathfrak{m}^{k}/\mathfrak{m}^{k+1}$.  
As is well-known, this ring is isomorphic to the polynomial ring 
$R=\Z[x_1,\dots,x_n]$.  Moreover, the Magnus embedding induces a graded
ring isomorphism $\gr \L \xrightarrow{\sim} \gr P=R$.  
Let $\hB $ be the $I$-adic completion of $B$, and  
$\gr \hB =\bigoplus _{k\ge 0} \mathfrak{m}^{k} \hB /\mathfrak{m}^{k+1}\hB $
the associated graded module.  Then, the canonical map $B\to \hB $ induces 
an isomorphism $\gr B \xrightarrow{\sim}\gr \hB $ of graded modules over 
the ring $R$.  

Combining these facts, we can restate Massey's result as follows:

\begin{thm}[{\cite{Ma}}]\label{thm:massey} The generating series for
the ranks of the Chen groups of $G$,
$\sum _{k=0}^{\infty }\theta _{k+2}t^{k},$
is equal to the Hilbert series of the graded module associated
to the $I$-adic completion of $B(G)$,
$\sum _{k=0}^{\infty }\operatorname{rank}
  (\mathfrak{m}^{k} \hB /\mathfrak{m}^{k+1}\hB )t^{k}.$
\end{thm}

An immediate consequence of this theorem is that, for $k$ sufficiently
large, $\theta _{k}$ is given by a polynomial in $k$.  Indeed, this is just 
the Hilbert-Serre polynomial of $\gr \hB $, see \cite{ZS}.

\subsection{Groebner Bases}\label{subsec:Groebner}
Let 
$\L ^{a} 
\xrightarrow{\Delta} \L ^{b} 
\rightarrow B 
\rightarrow 0$
be a (finite) presentation of the Alexander invariant.  
Note that, by replacing the generators of the free module $\L ^{a}$ by
suitable multiples if necessary, we may assume that the entries of
the matrix of $\Delta $  are polynomials in the variables $t_{i}$.
Let $J=\im \Delta $.  A presentation for the $I$-adic
completion of $B$ is given by
$\DS {\hL ^{a} \xrightarrow{\hD}\hL ^{b}
\rightarrow \hB \rightarrow 0}$,
where $\hD $ is obtained from $\Delta $ via the Magnus embedding.
Clearly, $\im \hD =\J $.  Since all the entries
of the matrix for $\hD $ belong to the subring $R\subset P$,
we may restrict $\hD $ to a map $\boldsymbol \Delta : R^{a}\to R^{b}$,
whose image, $\J \cap R ^{b}$, we denote by $\mathbf{J}$.

We must find a presentation for the associated graded module
$\gr \hB =\gr (P^{b}/\J )$.  This module is isomorphic to 
$R^{b}/\LT (\mathbf{J} )$, where $\LT (\mathbf{J} )$ is the 
submodule of $R^{b}$ consisting of lowest degree homogeneous 
forms of elements in $\mathbf{J}$, see \cite{ZS}.  We are left 
with finding a finite generating set for $\LT (\mathbf{J})$.
Such a set is provided by Mora's algorithm for
obtaining the tangent cone of an affine variety at the origin,
see \cite{CLO}, \cite{BW}.  Essentially, we must 
determine a (minimal) {\em Groebner basis}
$\G =\{g_{1},\dots ,g_{c}\}$ for the module $\mathbf{J}$, with
respect to a suitable monomial ordering.  Then, $\LT (\mathbf{J})$
has Groebner basis $\LT (\G )=\{\LT g_{1},\dots ,\LT g_{c}\}$, 
from which we can extract a minimal Groebner basis
$\H =\{h_{1},\dots ,h_{d}\}$.  
Putting all these facts together, 
we obtain the following.

\begin{thm}\label{thm:grB} The module $\gr \hB $ has presentation
$R^{d} 
\xrightarrow{\gr \hD} R^{b} 
\rightarrow \gr \hB 
\rightarrow 0$, where the rows of $\gr \hD $
constitute a minimal Groebner basis for the module generated by
the rows of the matrix $\hD $, obtained from a presentation
matrix $\D $ for $B$ by replacing $t_{i}$ by $1-x_{i}$.
\end{thm}

\begin{exm}\label{exm:Fn} Let $G = F_n $ be a finitely generated free 
group.  A presentation for the Alexander invariant $B$ of $G$ is given by 
$\L^a \xrightarrow{d_3} \L^b \to B \to 0$, where $a=\binom{n}{3}$, 
$b=\binom{n}{2}$, and $d_3$ is the differential in the standard 
$\L$-resolution of $\Z$.  In this instance, it is readily checked 
that the rows of the matrix $\hat d_3$ form a Groebner basis for the 
module $\mathbf{J}$.  A standard argument then yields the ranks of 
the Chen groups of $F_n$: $\theta_1 = n$ and $\theta_k=(k-1)\binom{k+n-2}{k}$ 
for $k\ge 2$, a calculation originally due to Murasugi \cite{Mu}.
\end{exm}

\subsection{Products}\label{subsec:Products}
Let $M_{1}$ and $M_{2}$ be two path
connected finite CW-complexes, with $K_{i}=H_{1}(M_{i})$
free abelian, and let $M'_{i}$ be the corresponding maximal
abelian covers.  Then $M=M_{1}\times M_{2}$ has maximal abelian
cover $\M =M'_{1}\times M'_{2}$, corresponding to
$K=H_{1}(M)=K_{1} \times K_{2}$.

\begin{prop}\label{prop:kunneth} There is an isomorphism of $\Z K$-modules,
\begin{equation*}  \label{eq:kunneth}
H_{1}(\M )\cong 
\bigl ( (H_{1}(M'_{1})\otimes _{\Z K_{1}} \Z K) \otimes _{\Z K_{2}} 
\Z \bigr ) \oplus \bigl ( (H_{1}(M'_{2})
\otimes _{\Z K_{2}} \Z K) \otimes _{\Z K_{1}} \Z \bigr ) .
\end{equation*}
\end{prop}
\begin{proof} 
By the K\"{u}nneth formula, the group $H_{1}(\M )$ is isomorphic to
$H_{1}(M'_{1})\otimes H_{0}(M'_{2}) 
  \oplus H_{0}(M'_{1})\otimes H_{1}(M'_{2}) .$
When viewed as a $\Z K$-module, the first summand is isomorphic to
\begin{equation*}  \label{eq:summand}
\begin{split}
&( (H_{1}(M'_{1})\otimes _{\Z K_{1}} \Z K)\otimes _{\Z K} \Z K_{1} )
  \otimes ( (\Z \otimes _{\Z K_{2}} \Z K)\otimes _{\Z K} \Z K_{2} ) \\
&\qquad = (H_{1}(M'_{1})\otimes _{\Z K_{1}} \Z K)  \otimes _{\Z K}
   (\Z K_{1}  \otimes \Z K_{2}) \otimes _{\Z K}
   (\Z \otimes _{\Z K_{2}} \Z K) \\
&\qquad = (H_{1}(M'_{1})\otimes _{\Z K_{1}} \Z K) \otimes _{\Z K_{2}} \Z ,
\end{split}
\end{equation*}
where we made use of the obvious isomorphism
$(\Z K_{1}  \otimes \Z K_{2})\cong \Z (K_{1}\times K_{2})$,
and viewed the induced module $H_{1}(M'_{1})\otimes _{\Z K_{1}} \Z K$ as
a $\Z K_{2}$-module by restriction of scalars.  The second summand
is treated exactly the same way.  
\end{proof}

We want to find now a presentation for the Alexander invariant
$B(M)=H_{1}(M' )$, given presentations for the Alexander invariants
$B(M_{i})=H_{1}(M'_{i})$.  Fix generators
$t_{1}^{(i)},\dots ,t_{n_{i}}^{(i)}$
for $K_{i}$, and use them to identify $\Z K_{i}$ with $\L _{i}$.

\begin{thm}\label{thm:alexprod}  If the Alexander invariants of $M_{1}$ 
and $M_{2}$ have presentations
$\L _{i}^{a_{i}} 
\xrightarrow{\D _{i}} \L _{i}^{b_{i}}
\rightarrow B(M_{i})
\rightarrow 0$,
then the Alexander invariant of $M=M_{1}\times M_{2}$ has presentation
\begin{equation*}  \label{eq:alexprod}
\L ^{a} 
\xrightarrow{\begin{pmatrix}\D _{1}\\
                    D_{2}^{b_{1}}
              \end{pmatrix}
              \oplus 
              \begin{pmatrix}\D _{2}\\ 
                     D_{1}^{b_{2}}\end{pmatrix}}
\L ^{b}
\rightarrow B(M)
\rightarrow 0,
\end{equation*}
where $a=a_{1}+n_{2}b_{1}+a_{2}+n_{1}b_{2},\, b=b_{1}+b_{2}$, and
$D_{i}=(t_{1}^{(i)}-1,\dots ,t_{n_{i}}^{(i)}-1)^{\top }$.
\end{thm}
\begin{proof}  Let us look at the first summand in the direct sum 
decomposition of $H_{1}(\M )$ from Proposition~\ref{prop:kunneth}.  
It is the tensor product over $\Z K$ of two
induced modules.  The first one is the $\Z K$-module induced from the
$\Z K_{1}$-module $H_{1}(M'_{1})$, and has presentation
\begin{equation}\label{eq:firstmod}
(\Z K)^{a_{1}} 
\xrightarrow{\D _{1}} (\Z K)^{b_{1}}
\rightarrow H_{1}(X_{1})\otimes _{\Z K_{1}}
\Z K
\to 0.
\end{equation}
The second one is the $\Z K$-module induced from the trivial
$\Z K_{2}$-module $\Z $, and has presentation
\begin{equation}\label{eq:secmod}
(\Z K)^{n_{2}} 
\xrightarrow{D_{2}} \Z K
\rightarrow \Z \otimes _{\Z K_{2}} \Z K
\rightarrow 0.
\end{equation}
Taking the tensor product (over $\Z K$) of the complexes \eqref{eq:firstmod}
and \eqref{eq:secmod} and truncating yields the following presentation
for the first summand of $H_{1}(\M )$:
\begin{equation*}\label{eq:firstsum}
(\Z K)^{a_{1}+n_{2}b_{1}} 
\xrightarrow{\begin{pmatrix}\D _{1}\\D_{2}^{b_{1}}\end{pmatrix}}
(\Z K)^{b_{1}}
\rightarrow (H_{1}(X_{1})\otimes _{\Z K_{1}} \Z K) \otimes _{\Z K_{2}} \Z 
\rightarrow 0.
\end{equation*}
The second summand is handled the same way, and that finishes the proof.
\end{proof}

\begin{cor}\label{cor:chenprod} 
The ranks of the Chen groups of $G=\pi_1(M_1 \times M_2)$ 
are given by
\begin{equation*}  \label{eq:chenprod}
\theta_k(G) = \theta_k(G_1)+\theta_k(G_2),
\end{equation*}
where $G_i=\pi_1(M_i)$.
\end{cor}

\begin{exm}\label{exm:prodFn} Let 
$G=F_{d_{1}}\times \dots \times F_{d_{\ell }}$ be a
direct product of finitely generated free groups.  Using the above
result, and the  calculation in Example~\ref{exm:Fn}, one can easily 
recover the ranks of the Chen groups of $G$ announced in \cite{CS1}:
$\theta _{1}=\sum _{i=1}^{\ell }d_{i}$ and
$\theta _{k}=(k-1)\cdot \sum _{i=1}^{\ell }\binom {k+d_{i}-2}k$ 
for $k\ge 2$.
\end{exm}

\begin{exm}\label{exm:cone} Let $\A$ be an affine arrangement of 
$n$ hyperplanes in $\C^d$, and let $\overline \A$ be the {\em cone} 
of $\A$, a central arrangement of $n+1$ hyperplanes in $\C^{d+1}$ 
(see \cite{OT}).  It is well-known that the complement $\overline M$ 
of $\overline \A$ is homeomorphic to the product of the complement 
$M$ of $\A$ and $\C^*$, $\overline M = M \times\C^*$.

Fix a generator $x$ for $\pi_1(\C^*) = \Z$.  Let 
$\L=\Z[t_1^{\pm 1},\dots,t_n^{\pm 1}]$, and suppose that the
Alexander invariant $B$ of $\A$ has presentation 
$\L^a \xrightarrow{\Delta} \L^b \to B \to 0$ (see 
section~\ref{sec:AlexInvArr}).  Using Theorem~\ref{thm:alexprod}, 
we obtain a presentation 
$\overline \L^{a+b} \xrightarrow{\overline \Delta} 
\overline \L^b \to \overline B \to 0$
for the Alexander invariant of $\overline \A$, where 
$\overline \L = \L[x^{\pm 1}]$ and 
$\overline \Delta = \begin{pmatrix}\D & (x-1)\cdot
\id\end{pmatrix}^{\top}$.   Thus the ranks of the Chen groups of
$G(\overline \A)$  coincide with those of $G(\A)$ for $k \ge 2$.  
Note that $\theta_1(\overline \A) = \theta_1(\A) + 1 = n+1$.

\end{exm}

\section{A Quick Trip through Fox Calculus}
\label{sec:FoxCalc}

In this section we review the basics of Fox's free differential calculus, 
as introduced in \cite{Fo}, and developed in \cite{Bi}, and derive some
consequences.    

\subsection{Fox Gradient}\label{subsec:FoxGrad} 
Let $F_{n}$ be the free group on generators $t_{1},\dots ,t_{n}$, 
and $\Z F_{n}$ its group ring.  Let $W_{n}=\bigvee _{1}^{n} S^{1}$ be
a wedge of $n$ circles, with basepoint $*$ at the wedge point.  
Let $\widetilde W_{n}$ be the universal cover, with basepoint 
$\widetilde *$, and let $\widetilde C_{\bullet }(\widetilde W_{n})$ 
be the augmented, equivariant chain complex of $\widetilde W_{n}$.  
Identifying $C_{0}(\widetilde W_{n})$ with $\Z F_{n}$, and 
$C_{1}(\widetilde W_{n})$ with $(\Z F_{n})^{n}$ (with basis 
$e_{1},\dots ,e_{n}$ given by the lifts of the 1-cells at $\widetilde *$),
we obtain the standard free $\Z F_{n}$-resolution of $\Z $, 
\begin{equation*} \label{eq:fnres}
0\rightarrow (\Z F_{n})^{n} 
\xrightarrow{\p_{1}} 
\Z F_{n} 
\xrightarrow{\epsilon}
\Z \rightarrow 0, 
\end{equation*}
where $\p _{1}(e_{i})=t_{i}-1$ and $\epsilon (t_{i})=1$.
The Fox Calculus is based on the observation that the augmentation ideal,
$IF_{n}=\ker \epsilon $, is a free $\Z F_{n}$-module of rank $n$,
generated by the entries of the matrix of $\p _{1}$.  This can be rephrased
as follows:  Given any $w\in \Z F_{n}$,  there exist unique elements
$\frac{\partial w}{\partial t_{i}} \in \Z F_{n}$
(called the Fox derivatives of $w$) such that the following ``\FTC '' holds:
\begin{equation}\label{eq:FTC}
w-\epsilon (w)=\sum _{i=1}^{n} \frac{\partial w}
{\partial t_{i}} (t_{i}-1).
\end{equation}

Let us define the {\em Fox gradient} to be the $\Z $-linear homomorphism
$\nabla :\Z F_{n} \to (\Z F_{n})^{n}$ given by
\begin{equation*}\label{eq:nabla}
\nabla (w) = \sum _{i=1}^{n} \frac{\p w}{\p t_{i}}\, e_{i}.
\end{equation*}
Then, formula \eqref{eq:FTC} takes the form
$\partial _{1}(\nabla (w)) = w - \epsilon (w).$
From this can be deduced the following ``product rule''
for the Fox gradient:
$\nabla (uv) = \nabla (u)\cdot \epsilon (v)+u\cdot \nabla (v).$
In particular, $\nabla (z^{-1})=-z^{-1}\nabla (z)$,
for $z\in F_{n}$.

Now consider an endomorphism $\alpha :F_{n}\to F_{n}$.  This defines
a map $\a :W_{n}\to W_{n}$ (unique up to homotopy).  The induced chain map
$\alpha _{\bullet }:C_{\bullet }
(\widetilde W_{n})\to C_{\bullet }(\widetilde W_{n})$
can be written as
\begin{equation}\label{eq:jack}
\begin{CD}
C_{1}(\widetilde W_{n})       @>\p _{1}>>             C_{0}(\widetilde W_{n})\\
@VVJ(\alpha )\circ \tilde \alpha V   @VV\tilde \alpha V \\
C_{1}(\widetilde W_{n})       @>\p _{1}>>             C_{0}(\widetilde W_{n})
\end{CD}
\end{equation}
where 
$J(\alpha ):(\Z F_{n})^{n}\to (\Z F_{n})^{n}$ is the {\em Fox Jacobian} of
$\a$; namely, the $\Z F_{n}$-linear homomorphism given by
$J(\a)(e_{i}) = \nabla (\a (t_{i})).$
If $\b :F_{n}\to F_{n}$ is another endomorphism, the fact that
$(\b\circ\a)_{\bullet }=\b _{\bullet }\circ\a _{\bullet}$
may be rephrased as the ``chain rule of Fox Calculus:''
$J(\a\cdot\b)=\tilde \b (J(\a))\cdot J(\b).$  
In particular, $J(\a)^{-1}=\tilde\a\circ J(\a^{-1})\circ\tilde\a^{-1}$. 

\subsection{Abelianized Fox Jacobian}
\label{AbelJack}
Let $\Z ^{n}$ be the free abelian group on generators
$t_{1},\dots ,\linebreak[0] t_{n} $, and identify the group ring
$\Z \Z ^{n}$ with $\L =\Z [t_{1}^{\pm 1},\dots ,t_{n}^{\pm 1}]$. Let 
$T^{n}=\times _{1}^{n} S^{1}$ be the $n$-torus.  
The augmented, equivariant chain complex,
$\widetilde C_{\bullet }=\widetilde C_{\bullet }(\widetilde T^{n})$,
of the universal (abelian) cover can be written as
\begin{equation}\label{eq:znres}
0\rightarrow C_{n} \xrightarrow{d_{n}} 
\cdots \rightarrow C_{3} \xrightarrow{d_{3}}
C_{2} \xrightarrow{d_{2}} 
C_{1} \xrightarrow{d_{1}} 
C_{0} \xrightarrow{\epsilon}
\Z  \rightarrow  0.
\end{equation}
Identifying $C_{0}$ with $\L $, $C_{1}$ with $\L ^{n}$, 
and $C_{k}$ with  $\bigwedge ^{k} C_{1}=\L^{\binom {n}k}$,
we obtain the standard free $\L$-resolution
of $\Z $, with differentials given by
$d_{k}(e_J)=\sum _{r=1}^{k} (-1)^{k+r} (t_{j_{r}}-1) 
\cdot e_{J\setminus\{j_r\}}$, where 
$e_J=e_{j_1}\wedge\dots\wedge e_{j_k}$ if
$J=\{j_1,\dots,j_k\}$. 

Let $\ab :F_{n}\to \Z ^{n}$, $x\mapsto x^{\ab }$, be the 
abelianization homomorphism.  
For an element $w\in \Z F_{n}$, let
$\nabla ^{\ab }(w):=\wab(\nabla (w))\in \L ^{n}$ be its
abelianized Fox gradient.  This defines a $\L $-linear homomorphism
\begin{equation*}\label{eq:nab}
\nabla ^{\ab }(w):C_{0}\to C_{1}, 
\qquad u\mapsto \nabla ^{\ab }(w)\cdot u.
\end{equation*}

For an endomorphism $\a$ of $F_{n}$, let
$\Theta (\a ):=\wab(J(\alpha )):C_{1}\to C_{1}$ be its abelianized
Fox Jacobian.  This is a $\L $-linear map, whose matrix
has rows $\Theta (\a )(e_{i})=\nabla ^{\ab }(\alpha (t_{i}))$.  
Abelianizing diagram~\eqref{eq:jack} yields the chain map:
\begin{equation}\label{eq:Theta}
\begin{CD}
C_{1}       @>d_{1}>>             C_{0}\\
@VV\Theta (\a )\circ\tilde\a V   @VV\tilde\a V \\
C_{1}       @>d_{1}>>             C_{0}
\end{CD}
\end{equation}

Set $\Theta _{k}(\a ) = \bigwedge ^{k} \Theta (\a ):C_k \to C_k$
(in particular, $\Theta_0=\id$).  A computation in the exterior 
algebra $C_{\bullet}=\bigwedge^{\bullet} C_1$ 
shows that $\Theta _{k-1}(\a )\circ\tilde\a \circ d_k = d_k \circ
\Theta _{k}(\a )\circ\tilde\a$ for each $k$, $1\le k \le n$.  Thus, 
\eqref{eq:Theta} extends to a chain map
\begin{equation}\label{eq:chainTheta}
\Theta_{\bullet}(\a)\circ \tilde\alpha : C_{\bullet}\to  C_{\bullet}
\end{equation}
This chain map is the composite of two chain maps.     
The first is the (non-$\L$-linear) map 
$\tilde\a: (C_{\bullet},d_{\bullet}) \to 
 (C_{\bullet},\tilde\a\circ d_{\bullet} \circ\tilde\a^{-1})$. 
The second is the ($\L$-linear) map $\Theta_{\bullet}(\a): 
(C_{\bullet},\tilde\a\circ d_{\bullet} \circ\tilde\a^{-1})\to
(C_{\bullet},d_{\bullet})$.

\subsection{Magnus representations}
\label{subsec:MagnusRep}
An automorphism $\a \in \Aut(F_{n})$ is called an 
{\em IA-automorphism} if its abelianization, 
$\ab (\a ):\Z ^{n}\to \Z ^{n}$,  is the identity map.  
In this case, $\tilde\a=\id$, and so 
$\Theta _{\bullet }(\a ):C_{\bullet }\to C_{\bullet }$ 
is a chain map.   

The set of IA-automorphisms forms a subgroup of 
$\Aut (F_{n})$, denoted by $\IA (F_{n})$. By the chain rule, 
$\Theta (\a \cdot \b) =\Theta (\a)\cdot \Theta (\b)$, 
for $\a,\b \in \IA (F_{n})$.  Thus 
$\Theta : \IA (F_{n}) \to \Aut _{\L }(C_{1})\cong \GL (n,\L )$ is
a linear representation of $\IA (F_{n})$, called the 
{\em Magnus representation}, see \cite{Bi}.  From the above 
discussion, we see that this representation generalizes to 
$\Theta _{k}: \IA (F_{n}) \to \Aut _{\L }(C_{k})
\cong \GL (\binom {n}{k},\L )$.  

\begin{rem}
For $\a \in \IA (F_{n})$, the chain automorphism  
$\Theta _{\bullet }(\a ):C_{\bullet }\to C_{\bullet }$ 
admits the following topological interpretation.  
The map $\a : W_{n} \to W_{n}$ lifts to a map of the maximal 
abelian covers, $\a ':W_{n}' \to W_{n}'$.  View 
$W'_{n}$ as the $1$-skeleton of $\widetilde T^{n}$.  
The map $\a '$ extends to a $\Z^n$-equivariant map 
$\bar \a : \widetilde T^{n} \to \widetilde T^{n}$.  
The induced chain map,  $\bar\a_{\bullet}: 
C_{\bullet }(\widetilde T^{n})\to C_{\bullet }(\widetilde T^{n})$, 
is chain-equivalent to $\Theta _{\bullet }(\a )$. 
\end{rem}

\section{The Alexander Invariant of a Free Automorphism}
\label{sec:AlexFreeAuto}

In this section, we find presentations for the Alexander module 
and the Alexander invariant of the group of an IA-automorphism.  
A more explicit presentation for the latter
is given in case the automorphism 
is basis-conjugating.  

\subsection{The Group of a Free Automorphism}
\label{subsec:GroupFreeAuto} 
Associated to an automorphism $\a$ of the free group
$F_{n}=\langle t_{1},\dots ,t_{n}\rangle $ is the group
\begin{equation*}\label{eq:Galpha}
G(\a )=\langle t_{1},\dots ,t_{n} \mid t_{1}=\a (t_{1}),\dots ,
t_{n}=\a (t_{n})\rangle .
\end{equation*}
Notice that $\a$ induces the identity automorphism on $G(\a)$.
In fact, $G(\a)$ is the maximal quotient of $F_{n}$ with
this property.  Also, note 
that $G(\a )$ is independent of
the choice of free generators for $F_{n}$:  If $x_{1},\dots , x_{n}$
is another such choice, then
$\langle t_{1},\dots ,t_{n} \mid t_{i}=\a (t_{i})\rangle 
\cong \langle x_{1},\dots ,x_{n} \mid x_{i}=\a (x_{i})\rangle$.
Finally, notice that the group of a free automorphism depends only on
the conjugacy class of that automorphism:  If $\b \in \Aut (F_{n})$, then
$G(\b^{-1} \circ \a \circ \b ) = \langle t_{i} \mid t_{i}=
\b^{-1} \circ \a \circ \b (t_{i})\rangle =
\langle t_{i} \mid \b (t_{i})=\a (\b (t_{i}))\rangle 
=\langle x_{i} \mid x_{i}=\a (x_{i})\rangle \cong G(\a )$.
See~\cite{Mo} for details.  

Topologically, the group $G(\a )$ can be interpreted as follows.
Recall that $W_{n}$ denotes a wedge of $n$ circles, and that
$\a :W_{n}\to W_{n}$
also denotes a basepoint preserving homotopy equivalence
that induces $\a :F_{n} \to F_{n}$ on fundamental groups.  Let
$Y(\a )=W_{n} \times _{\a } S^{1}$ be the mapping torus of $\a $; its
fundamental group is the semidirect product
$F_{n}\rtimes _{\a } \Z =\langle t_{1},\dots , t_{n}, x\mid x^{-1} t_{i} x =
\a (t_{i})\rangle $. Let
$X(\a ) = W_{n} \times _{\a } S^{1} \bigcup _{*\times S^{1}} *\times D^{2}$.
Then $\pi _{1}(X(\a ))=G(\a )$, and, in fact, $X(\a )$ is homotopy
equivalent to the 2-complex associated to the above
presentation of $G(\a )$.

\subsection{Alexander Invariants}
\label{subsec:AlexInvIA}
Let $\a $ be an $\IA $-automorphism 
of $F_{n}$, and $G=G(\a )$ the associated group.  Then $H_{1}(G)=\Z ^{n}$, 
the free abelian group generated by $t_{1},\dots , t_{n}$.
Let $p:\X \to X$ be the corresponding (maximal abelian) 
cover of $X=X(\a )$.
We call the $\L $-modules $A(\a )=H_{1}(\X,p^{-1}(*) )$, resp. 
$B(\a )=H_{1}(\X )$ the {\em Alexander module}, resp. 
{\em Alexander invariant} of $G(\a )$.    
We wish to find presentations for these modules.  

First consider $Y=W_{n} \times _{\a } S^{1}$.  The chain complex of 
its maximal abelian cover is obtained using the Fox calculus
as in \cite{CS2}:
\begin{equation*}\label{eq:chainY}
C_{\bullet }(Y'):\qquad C_{2}(Y') 
\xrightarrow{\begin{pmatrix}\id - x\cdot \Theta (\a ) & d_{1} \end{pmatrix}}
C_{1}(Y') 
\xrightarrow{\begin{pmatrix}d_{1} \\ x-1 \end{pmatrix}}
C_{0}(Y') 
\xrightarrow{\epsilon}\Z 
\rightarrow 0
\end{equation*}
where the chain groups are the modules over $\overline{\L}=\L [x^{\pm 1}]$
given by $C_{2}(Y')=C_{1}\otimes _{\L }\overline{\L}$,\
$C_{1}(Y')=(C_{1}\oplus C_{0})\otimes _{\L }\overline{\L}$, and
$C_{0}(Y')=C_{0}\otimes _{\L }\overline{\L}$.
It follows that the chain complex of the maximal abelian 
cover of $X=X(\a )$ is
\begin{equation}\label{eq:alexmod}
C_{\bullet }(\X ):\quad C_{1} 
\xrightarrow{\id -\Theta (\a )} C_{1}
\xrightarrow{d_{1}} C_{0} 
\xrightarrow{\epsilon} \Z 
\rightarrow 0.
\end{equation}
Hence, $A(\a )=\coker (\id -\Theta (\a ))$ and
$B(\a )=\ker (d_{1})/\im (\id -\Theta (\a )).$
By homological algebra, there exists a chain map from the 
chain complex \eqref{eq:alexmod} to the free $\L $-resolution 
\eqref{eq:znres}, extending the identity map of $\Z $:
\begin{equation*}\label{eq:CDPhi}
\begin{CD}
@. C_{1} @>\id -\Theta (\a )>> C_{1} @>d_{1}>>
C_{0} @>\epsilon >>\Z @>>> 0\\
@. @VV{\Phi (\a )}V @VV{=}V @VV{=}V @VV{=}V \\
C_{3} @>d_{3}>> C_{2} @>d_{2}>> C_{1} @>d_{1}>> C_{0}
 @>\epsilon >>\Z @>>> 0\
\end{CD}
\end{equation*}
A diagram chase shows that
$B(\a )=\coker  \bigl(\begin{smallmatrix}\Phi (\a ) \\ 
d_{3}\end{smallmatrix}\bigr).$  
To summarize, we have:

\begin{prop}\label{prop:IAB}  If $\a \in \IA (F_{n})$,
the Alexander invariants of $G(\a )$ have presentation
\begin{equation*}  \label{eq:IAB}
C_{1}
\xrightarrow{\id - \Theta(\a)}
C_{1}
\rightarrow{A(\a )}
\rightarrow 0, 
\quad\text{resp.}\quad
C_{1}\oplus C_{3} 
\xrightarrow{\begin{pmatrix}\Phi (\a ) & d_{3}\end{pmatrix}^\top}
C_{2}
\rightarrow{B(\a )}
\rightarrow 0.
\end{equation*}
\end{prop}

\begin{rem}\label{rem:modd3}  The map $\Phi (\a )$ is not unique, 
but rather, it is unique up to chain homotopy:  Given two choices,
$\Phi _{1}(\a )$ and $\Phi _{2}(\a )$, there is a homomorphism
$D:C_{1}\to C_{3}$ such that $\Phi _{1}(\a )-\Phi _{2}(\a )=
d_{3}\circ D$.  We abbreviate this by saying that
$\Phi_{1}(\a)=\Phi_{2}(\a)\mod d_{3}$.  Of course, any
two choices yield equivalent presentations for $B(\a)$.
\end{rem}

As noted previously, the group $G(\b \circ \a \circ \b ^{-1})$ is
isomorphic to $G(\a )$.  The relationship between the
corresponding chain maps is as follows:
\begin{prop}\label{prop:Phiconj}  For $\a , \b \in \IA (F_{n})$,
$
\Phi (\b \circ \a \circ \b ^{-1})=
\Theta _{2}(\b )\circ \Phi (\a )\circ 
\Theta (\b ^{-1}) \mod d_{3}.
$
\end{prop}
\begin{proof} By Remark~\ref{rem:modd3}, it is enough to show that
\begin{equation*}\label{eq:bab'}
d_{2}\circ \Theta _{2}(\b )\circ \Phi (\a )\circ \Theta (\b ^{-1}) =
\id -\Theta (\b \circ \a \circ \b ^{-1}).
\end{equation*}
Since the right-hand side equals
$\Theta (\b )\circ (\id -\Theta (\a ))\circ \Theta (\b ^{-1})$,
the claim follows from the equalities 
$d_{2}\circ \Theta _{2}(\b )=\Theta (\b )\circ d_{2}$
and $d_{2}\circ \Phi (\a )=\id -\Theta (\a )$.
\end{proof}

\subsection{Basis-Conjugating Automorphisms}\label{subsec:BasisConj}
An automorphism $\a $ of $F_{n}=\langle t_{1},\dots ,t_{n}\rangle $
is called a {\em basis-conjugating automorphism} if there exists an
$n$-tuple $\mathbf{z} = (z_{1}, \dots , z_{n})$, with $z_{i}\in F_{n}$, 
such that $\a =\c _{\mathbf{z}}$, where
$\c _{\mathbf{z}}(t_{i})=z_{i}t_{i}z_{i}^{-1}.$
The basis-conjugating automorphisms of $F_{n}$ form a subgroup,
$\CA (F_{n})$, 
of $\Aut (F_{n})$.  For $\a \in \CA (F_n)$, the
following definition/proposition gives an explicit
formula for $\Phi (\a )$.
\begin{prop}\label{prop:Phidef} For $\c _{\mathbf{z}}\in \CA (F_{n})$,
define the $\L $-homomorphism
$\Phi (\c _{\mathbf{z}}):C_{1}\to C_{2}$ by
\begin{equation}  \label{eq:Phidef}
\Phi (\c _{\mathbf{z}})(e_{i}) = \nabla ^{\ab }(z_{i})\wedge e_{i}.
\end{equation}
Then $\id - \Theta (\c _{\mathbf{z}})=d_{2}\circ \Phi (\c _{\mathbf{z}}).$
\end{prop}
\begin{proof}   First, note that the Magnus representation of
$\c _{\mathbf{z}}$ is given by:
\begin{equation*}\label{eq:magnus}
\Theta (\gamma _{\mathbf{z}})(e_{i})=(1-t_{i})\cdot \nabla ^{\ab }(z_{i}) +
z_{i}^{\ab }\cdot e_{i}.
\end{equation*}
Hence:
\begin{align*}
\left ( \id - \Theta (\gamma _{\mathbf{z}}) \right ) (e_{i})
&=(t_{i}-1)\cdot \nabla ^{\ab }(z_{i}) + (1-z_{i}^{\ab })\cdot e_{i}\\
&=d_{1}(e_{i})\cdot 
   \nabla ^{\ab }(z_{i}) - d_{1}(\nabla ^{\ab }(z_{i}))\cdot e_{i}\\
&=d_{2}(\Phi (\gamma _{\mathbf{z}})(e_{i})).   
\qquad\qquad \qed
\end{align*}
\renewcommand{\qed}{}\end{proof}

\begin{rem}  \label{rem:Bpurelink}
As mentioned before, an explicit formula for the Alexander 
invariant $B(L)$ of an arbitrary link $L\subset S^3$ is lacking.  
If $L$ is a pure link of $n$ components, though, 
Propositions~\ref{prop:IAB} and ~\ref{prop:Phidef} provide a 
presentation for $B(L)$, with $\binom{n}{2}$ generators 
and $\binom{n}{3}+n$ relations.  Indeed, 
as shown by Artin, the braid group $B_n$ admits a faithful 
representation $B_n \hookrightarrow \Aut(F_n)$, 
which restricts to $P_n\hookrightarrow  \CA(F_n)$.  
Moreover, any link $L$ is the closure, $\hat\alpha$, of a braid $\a\in B_n$, 
and $\pi_1(S^3\setminus L) = G(\a)$.  Now assume $L$ is a pure link,  
i.e,~$L=\hat\a$, for some $\a\in P_n$.  
Then $\a=\gamma_{\mathbf{z}}$, where $z_i$ is the longitude corresponding 
to the meridian $t_i$, and we get 
$B(L)=\coker  \bigl(\begin{smallmatrix}\Phi (\gamma_{\mathbf{z}}) \\ 
d_{3}\end{smallmatrix}\bigr)$, with $\Phi (\gamma_{\mathbf{z}})$ 
given by~\eqref{eq:Phidef}. 
\end{rem}

\subsection{Alexander Invariant of Several Automorphisms}
\label{subsec:Generalization} 
The above notions generalize in a straightforward manner, from a single
automorphism $\a $ to several automorphisms $\a _{1},\dots ,\a _{s}$ of
$F_{n}$.  Namely, let
\begin{equation*}\label{eq:Gmono}
G(\a _{1},\dots ,\a _{s})=\langle t_{1},\dots , t_{n}
\mid t_{i}=\a _{k}(t_{i}),\; 1\le k\le s\rangle 
\end{equation*}
be the maximal quotient of $F_{n}$ on which all $\a _{k}$ act trivially.
This group can also be characterized as the quotient of
the semidirect product $F_{n} \rtimes F_{s} = \langle t_{1},\dots ,t_{n},
x_{1},\dots ,x_{s} \mid x_{k}^{-1}t_{i}x_{k}
=\a _{k}(t_{i})\rangle $ by the normal
closure of $F_{s}=\langle x_{1},\dots ,x_{s}\rangle $.

Assume $\a _{k}\in \IA (F_{n})$, for $1\le k\le s$.  Let $Y$ be the
presentation $2$-complex for $F_{n} \rtimes F_{s}$, and $X$ that of
$G=G(\a _{1},\dots ,\a _{s})$.  The chain complex of the maximal 
abelian cover of $Y$ has the following form:
\begin{equation*}  \label{eq:chainbigY}
C_{\bullet}(Y'):\quad C_{2}(Y') 
\xrightarrow{%
  \begin{pmatrix}\id -x_{1}\cdot \Theta (\a _{1}) & d_{1} &\cdots & 0 \\ 
  \vdots & & \ddots \\ 
   \id -x_{s}\cdot \Theta (\a _{s}) & 0 &\cdots & d_{1}
  \end{pmatrix}}
C_{1}(Y') 
\xrightarrow{\begin{pmatrix}d_{1} \\ \overline{d}_{1}\end{pmatrix}}
C_{0}(Y') 
\xrightarrow{\epsilon}
 \Z 
\rightarrow 0
\end{equation*}
where the chain groups are the modules over
$\overline{\L}=\L [x_{1}^{\pm 1},\dots ,x_{s}^{\pm 1}]$
given by
$C_{2}(Y')=\oplus _{1}^{s} C_{1}\otimes _{\L }\overline{\L}$,\
$C_{1}(Y')=(C_{1}\oplus \oplus _{1}^{s} C_{0})\otimes _{\L }\overline{\L}$, and
$C_{0}(Y')=C_{0} \otimes _{\L }\overline{\L}$,
and where 
$\overline{d}_{1}=\begin{pmatrix}x_{1}-1 & \cdots & x_{s}-1 
\end{pmatrix}^{\top }$.
The chain complex of the maximal abelian cover of $X$ is then
\begin{equation*}\label{eq:chainbigX}
C_{\bullet}(\X ):\quad C_{1}^{s} 
\xrightarrow{\begin{pmatrix}\id - \Theta (\a _{1})\\ \cdots \\
             \id - \Theta (\a _{s})\end{pmatrix}} 
C_{1} 
\xrightarrow{d_{1}} C_{0} 
\xrightarrow{\epsilon} \Z 
\rightarrow 0.
\end{equation*}

This chain complex provides a presentation 
matrix---the so-called {\em Alexander matrix}, 
$\left( \id-\Theta(\a _{1}) \ \cdots 
\ \id-\Theta(\a _{s}) \right) ^{\top }$---for 
the Alexander module $A=A(\a_1,\dots,\a_s)$.  
Furthermore, if $\Phi (\a _{k}):C_{1}\to C_{2}$ 
satisfy $d_{2}\circ \Phi (\a _{k})=\id - \Theta (\a _{k}) $, 
then the Alexander invariant $B=B(\a_1,\dots,\a_s)$ 
has presentation matrix 
$\begin{pmatrix}\Phi (\a _{1}) &\cdots &\Phi (\a _{s}) & d_{3}
\end{pmatrix}^{\top }$.   When $\a_k\in\CA(F_n)$, we obtain 
an explicit presentation for $B$.

\begin{thm}\label{thm:Bpres}   
Let $\c _{\mathbf{z}^{1}},\dots ,\c _{\mathbf{z}^{s}}$
be a collection of basis-conjugating automorphisms of $F_{n}$.  
Let $\Phi(\a_k):C_1 \to C_2$ be the homomorphisms defined 
by \eqref{eq:Phidef}.  Then the Alexander invariant of 
$G(\c _{\mathbf{z}^{1}},\dots ,\c _{\mathbf{z}^{s}})$
has presentation
\begin{equation*}   \label{eq:Bpres}
C_{1}^{s}\oplus C_{3} 
\xrightarrow{\begin{pmatrix}\Phi (\c _{\mathbf{z}^{1}}) &\cdots &
         \Phi (\c _{\mathbf{z}^{s}}) & d_{3}
         \end{pmatrix}^{\top } }
C_{2}
\rightarrow B
\rightarrow 0.
\end{equation*}
\end{thm}

\section{Local Alexander Invariants}
\label{sec:LocalAlexInv}

We now find presentations for the Alexander invariant of the group of a 
full-twist braid automorphism, and that of a related ``vertex" group.  
These presentations are given in terms of the Gassner representation, 
$\Theta:P_n\to \GL(n,\L)$, which is simply the restriction of the 
Magnus representation, $\Theta:\IA(F_n)\to \GL(n,\L)$, to the 
pure braid group $P_n$.

\subsection{Alexander Invariant of a Twist Automorphism}
\label{subsec:TwistAuto}
Let $V=\{i_{1},\dots ,i_{r}\}$ be an increasingly ordered subset of 
$[n]=\{1,\dots,n\}$.  Let $A_V$ the pure braid in $P_n$ which 
performs a full twist on the strands corresponding to $V$, leaving 
the other strands fixed.  Let $\sigma_i$ ($1\le i <n$) be the standard 
generators of $B_n$, and $A_{i,j}=\sigma _{j-1}\cdots \sigma _{i+1}
\sigma _{i}^{2}\sigma _{i+1}^{-1}\cdots \sigma _{j-1}^{-1}$ 
($1\le i<j\le n$) the standard generators of $P_n$, see~\cite{Bi}.  
The twist on $V$ is given by 
\begin{equation}\label{eq:A_V}
A_{V} = (A_{i_{1},i_{2}})(A_{i_{1},i_{3}}A_{i_{2},i_{3}})(A_{i_{1},i_{4}}
A_{i_{2},i_{4}}A_{i_{3},i_{4}})\cdots 
(A_{i_{1},i_{r}}\cdots A_{i_{r-1},i_{r}}). 
\end{equation}
A computation with the Artin 
representation reveals that $A_{V}=\c _{\mathbf{w}}$, where 
$\mathbf{w} = (w_{1},\dots ,w_{n})$ is defined as follows:
\begin{equation}\label{eq:w_i}
w_{i}=\begin{cases}t_{V}^{} &\text{ if } i\in V,\\ \relax
 [t_{V^{i}}^{},t_{^{i}V}^{}]
&\text{ if } i\in \overline{V}\setminus V,\\
 1 &\text{ otherwise,}
\end{cases}
\end{equation}
where $\overline{V} = \{i\in [n] \mid i_{1}\le i \le i_{r}\}$,
$V^{i}=\{ j\in V \mid j<i\}$,
$^{i}V=\{j\in V\mid i<j\}$, and
$t_{V}^{}=\prod _{j\in V} t_{j}=t_{i_1}\cdots t_{i_r}$.  

Let $G(A_V)=\langle t_1,\dots,t_n \mid A_V(t_i)=t_i\rangle$ be the group 
associated to $A_V\in\Aut(F_n)$.  A computation with~\eqref{eq:w_i} 
shows that 
\begin{equation*}\label{eq:GAV}
G(A_V)=\langle t_1,\dots,t_n \mid 
t_V^{} t_i t_V^{-1} = t_i\ \ , i\in V\rangle, 
\end{equation*}
and so $G(A_V)\cong (F_{r-1}\times F_1)*F_{n-r}$. 
By Proposition~\ref{prop:IAB}, the Alexander invariant of $G(A_V)$ 
has presentation
\begin{equation}  \label{eq:BAV}
C_{1}\oplus C_{3} 
\xrightarrow{\begin{pmatrix}\Phi(A_V) & d_{3}\end{pmatrix}^\top}
C_{2}
\rightarrow{B(A_V)}
\rightarrow 0.
\end{equation}
By Proposition~\ref{prop:Phidef}, and a Fox calculus computation, 
the map $\Phi _{V}^{}:=\Phi (A_{V}):C_{1} \to C_{2}$ is given by
\begin{equation}\label{eq:PhiV}
\Phi _{V}^{}(e_{i})
    =\begin{cases}\nabla _{V} \wedge e_{i} &\text{ if } i\in V\\
      (1-t_{i})\nabla _{V}\wedge \nabla _{^{i}V}
      &\text{ if } i\in \overline{V}\setminus V\\
      0 &\text{ otherwise.}
      \end{cases}
\end{equation}
where 
$\nabla _{V}:=\nabla ^{\ab }(t_{V})=\sum _{i\in V}t_{V^{i}}e_{i}$.

\subsection{Simplified Presentation for $B(A_V)$}
\label{subsec:BAV}  
Set $C_{k}(V)=\Span  \{ e_{J}\mid J\subset V\}$, and let 
$\iota_{V}:C_k(V)\to C_k$ be the inclusion, and 
$\pi_{V}:C_k\to C_k(V)$ the natural projection.  
Write $V':=V\setminus \{\min V\}=\{i_{2},\dots , i_{r}\}$.  
From~\eqref{eq:PhiV} it is apparent that 
$\Phi_V(C_1)=\Phi_V(C_1(V))\subset C_2(V)$.  
Since $0=\nabla_V \wedge \nabla_V = \nabla_V \wedge e_{i_1} + 
\nabla_V \wedge \sum _{i\in V'}t_{V^{i}}e_{i}$, we see that 
$\Phi_V(e_{i_1}) \in \Phi_V(C_1(V'))$.  Thus,
\begin{equation}  \label{eq:imPhi}  
\Phi _{V}(C_{1})=\Phi _{V}(C_{1}(V'))\subset C_{2}(V).
\end{equation}

Define an automorphism $\mu_V: F_n \to F_n$ by:
\begin{equation*}  \label{eq:muV}  
\mu_{V}(t_i)=
  \begin{cases}
     t_{V} &\text{ if } i=\min V,\\
     t_i &\text{ otherwise.}
  \end{cases}
\end{equation*}
Note that $\Theta(\mu_V)(e_{i_1})=\nabla_V$, and 
$\Theta_2(\mu_V)(e_{i_1}\wedge e_i) =\Phi_V(e_i)$, for $i\in V'$. 
Thus, 
\begin{equation}  \label{eq:muPhi}  
\Theta_2(\mu _{V})^{-1}\circ \Phi _{V}(e_{i}) = e_{i_{1}}\wedge e_{i} 
\quad \text{for } i\in V'.
\end{equation}

Let $C'_2(V)$ be the direct summand of $C_2$ spanned by 
$\{e_{i_{1}}\wedge e_{i} \mid i\in V'\}$, let $C_2^{\perp}(V)$ be the 
complementary summand,  and let $\pi^V:C_2 \to C_2^{\perp}(V)$ 
be the canonical projection.  Putting together \eqref{eq:BAV}, 
\eqref{eq:imPhi}, and \eqref{eq:muPhi}, we obtain:

\begin{prop}\label{prop:alexGAV}  
The  Alexander invariant of $G(A_V)$ has presentation
\begin{equation}  \label{eq:alexGAV}
C_{3} 
\xrightarrow{\Delta(V)}  
C_2^{\perp}(V)
\rightarrow B(A_V)
\rightarrow 0, 
\end{equation}
where $\Delta(V)=\pi^V \circ \Theta_2(\mu_V)^{-1} \circ d_{3}$.
\end{prop}

\subsection{Alexander Invariant of a Vertex Group}
\label{subsec:BV}  
To the twist automorphism $A_V$, we also associate a ``vertex group,"
$G_V:=G(\{A_V ,A_{i,j} \mid |\{i,j\}\cap V| \le 1\})$.  
Using \eqref{eq:w_i}, we obtain the following presentation: 
\begin{equation}  \label{eq:GV}  
G_V= \langle t_1,\dots ,t_n \mid t_V^{} t_i t_V^{-1} = t_i
\text{ if } i\in V, \ \
 t_j^{} t_i t_j^{-1} = t_i
\text{ if } \{i,j\}\not\subset V\rangle. 
\end{equation}

A (minimal) presentation for the Alexander invariant $B_V=B(G_V)$ 
may be obtained from~\eqref{eq:alexGAV} by restricting the map 
$\Delta(V)$ to a map $C_2(V') \wedge C_1 \to C_2(V')$, 
and some further matrix operations.  Alternatively, it may 
be obtained by applying Theorem~\ref{thm:alexprod} to 
the direct product decomposition $G_V\cong F_{r-1}\times\Z^{n-r+1}$,  
apparent from \eqref{eq:GV}.  The result is as follows:

\begin{prop} The  Alexander invariant of $G_V$ has presentation
\begin{equation*} \label{eq:presBV}
C_2(V')\wedge C_1
\xrightarrow{\Delta_V}
C_2(V')
\rightarrow
B_V
\rightarrow
0, 
\end{equation*}
where $\Delta_V=\pi_{V'}\circ \tilde\mu_V^{} \circ d_3 \circ 
\tilde\mu_V^{-1}\circ (\iota_{V'}\wedge \id)$.   
\end{prop}
\noindent
(This will be useful only when $|V|\ge 3$; if $|V|=2$, then 
$G_V=\Z^n$, $B_V=0$, and $C_2(V')=0$.)

The above presentation
may be extended to a free resolution, 
\begin{equation} \label{eq:BVres}
\dots
\rightarrow
C_2(V') \wedge C_{2}
\xrightarrow{\Delta^2_V}
C_2(V') \wedge C_{1}
\xrightarrow{\Delta_V}
C_2(V')
\rightarrow
B_V
\rightarrow
0,
\end{equation}
with boundary maps $\Delta^{\bullet}_V$ given by
\begin{equation*} \label{eq:DeltakV}
\Delta^k_V  = (\pi_{V'}\wedge \id)\circ\tilde\mu_V\circ d_{k+2}\circ 
\tilde\mu_V^{-1}\circ (\iota_{V'}\wedge \id) : 
C_2(V')\wedge C_{k} \to C_2(V')\wedge C_{k-1}.
\end{equation*}
Furthermore, by the discussion following \eqref{eq:chainTheta}, 
there exists a naturally defined chain map 
$\Psi_{V,\bullet}:(C_{\bullet}, d_{\bullet}) \to  
(C_2(V')\wedge C_{\bullet-2}, \Delta_V^{\bullet-2})$, given by 
\begin{equation} \label{eq:ChainV}
\Psi_{V,k}=(\pi_{V'}\wedge \id)\circ \Theta_k(\mu_V)^{-1}: 
C_k \to C_2(V')\wedge C_{k-2}, \ \ \text{ for } k\ge 2.
\end{equation}

\section{The Alexander Invariant of an Arrangement}
\label{sec:AlexInvArr}

In this section, we use the results of the previous sections
to obtain a presentation for the Alexander invariant of the
group of a hyperplane arrangement. 

\subsection{Braid Monodromy} \label{subsec:bmono}
The fundamental group of the complement of an arrangement of
complex hyperplanes is, by a Lefschetz-type theorem
of Zariski, isomorphic to that of a generic two-dimensional
section.  So, for the purpose of computing the Alexander invariant, it
is enough to consider affine line arrangements in $\C ^{2}$.
Let $\A =\{H_{1},\dots ,H_{n}\}$ be such an arrangement, 
with vertices 
$\mathcal V = \{v_{1},\dots ,v_{s}\}$.  
If $v_k = H_{i_1}\cap \dots \cap H_{i_r}$, let $V_k=\{i_1,\dots,i_r\}$ 
denote the corresponding ``vertex set.''  We identify the 
set $L_2(\A)$ of rank two elements in the lattice of $\A$ and the 
collection $\{V_1,\dots,V_s\}$ of vertex sets of $\A$.

The braid monodromy of $\A$ is determined as follows (see~\cite{CS3}
for details).  Choose coordinates $(x,z)$ in $\C ^{2}$ so that 
the projection $\pr _{1}:\C ^{2}\to \C $ is generic
with respect to $\A $.   Let $f(x,z)=\prod _{i=1}^{n} (z-a_{i}(x))$
be a defining polynomial for $\A $.  The root map
$a=(a_{1},\dots , a_{n}):\C \to \C ^{n}$ restricts to a map
from the complement of $\mathcal{Y}=\pr _{1}(\mathcal{V})$ to the
complement of the braid arrangement 
$\A _n=\{\ker (y_{i}-y_{j})\}_{1\le i <j \le n}$.
Identify $\pi _{1}(\C \setminus \mathcal{Y})$ with the free
group $F_{s}=\langle x_{1},\dots ,x_{s}\rangle $, and
$\pi _{1}(\C ^{n} \setminus \A _{n})$ with the pure braid group $P_{n}$.
Then, the {\em braid monodromy} of $\A $ is the induced
homomorphism on fundamental groups, $\a :F_{s} \to P_{n}$.

The braid monodromy generators $\a _{k}=\a (x_{k})$ can 
be written explicitly using a {\em braided wiring diagram} $\W$ 
associated to $\A$.    
Such a diagram, determined by the choices made above, may be (abstractly)
specified by a sequence of vertex sets and braids,
$\W =\W_s= \{ V_1,\b_1,V_2,\b_2,\dots,\b_{s-1},V_s\}$.
The braid monodromy generators are given by
$\a _{k}=A_{V_{k}}^{\d _{k}}$, where 
$A_{V_k}$ is the twist braid defined in~\eqref{eq:A_V} and 
$\d_k$ is a pure braid determined by the subdiagram $\W_k$.

\subsection{The Presentation for $B(\A)$}
\label{subsec:PresB}
Let $M=M(\A)$ be the complement of $\A$.
Let $G=G(\a_1,\dots , \a_s)$ be the fundamental group of $M$, 
with Alexander invariant $B=B(\A)$.
Theorem~\ref{thm:Bpres} 
provides the following presentation for $B$:
\begin{equation*}\label{eq:Barrangement}
C_{1}^{s}\oplus C_{3} 
\xrightarrow{\begin{pmatrix}\Phi _{V_{1}}^{\d _{1}}  &\cdots &
         \Phi _{V_{s}}^{\d _{s}} & d_{3}
   \end{pmatrix}^{\top }} 
C_{2}
\rightarrow B
\rightarrow 0, 
\end{equation*}
where 
$\Phi _{V}^{\d}:=\Phi (A_{V}^{\d })=\Theta _{2}(\d )\circ \Phi _{V}^{}
\circ \Theta (\d ^{-1}):C_{1} \to C_{2}$, 
and $\Phi_V$ is given by \eqref{eq:PhiV}.

This presentation can be simplified, based on the following
elementary observation:  If $R$ is a ring, and $B$ is an $R$-module,
with presentation 
$R^{p} \xrightarrow{\Delta} R^{q} 
\rightarrow B\rightarrow 0$, where
$\Delta =\Upsilon \circ \Xi $, or $\Delta =\Xi \circ \Upsilon $,
with $\Xi $ invertible, then $B$ can also be presented as
$R^{p} \
\xrightarrow{\Upsilon} R^{q} 
\rightarrow B
\rightarrow 0$.  
Since the maps $\Theta (\d _{k}^{-1})$ are invertible
we may replace $\Phi_{V_k}^{\d_k}$
by $\Theta_{2}(\d _{k})\circ \Phi _{V_{k}}$.  Furthermore, 
by~\eqref{eq:imPhi}, we may subsequently restrict each of the maps
$\Phi_{V_k}$ to $\Phi_k:C_1(V_k') \to C_2$.  Thus, we obtain
the following:

\begin{thm}\label{thm:alexarr}  The  Alexander invariant of an 
arrangement $\A $, with braid monodromy generators 
$A_{V_{1}}^{\d _{1}},\dots, A_{V_{s}}^{\d _{s}}$, has presentation
\begin{equation*}  \label{eq:alexinv}
K_1
\xrightarrow{\Delta} 
K_0
\rightarrow B(\A )
\rightarrow 0, 
\end{equation*}
where $\displaystyle{K_1=\bigoplus_{k=1}^s C_{1}(V_k') \oplus C_{3}}$, 
$K_0=C_2$, 
and $\Delta = \begin{pmatrix} \Phi & d_3 \end{pmatrix}^{\top}$, 
with $\Phi|_{C_1(V_k')}=\Theta _{2}(\d _{k})\circ \Phi_k$.
\end{thm}
Note that this presentation has $\binom {n}{2}$ generators and 
$\sum _{k=1}^{s}(|V_{k}|-1)+\binom {n}{3}$ relations, and that
$\sum _{k=1}^{s}(|V_{k}|-1)=b_{2}(M)$.  

\subsection{Real Arrangements}
\label{subsec:RealArrangements}
The presentation can be simplified in the case where $\A $ is 
the complexification of a line arrangement $\A_\R$ in $\R^2$.  
In this instance, the wiring diagram $\W$ can be chosen so that 
it contains no intermediary braids, and each ``conjugating braid,'' $\d_k$, 
is a subword of the full twist, $A_{[n]}$, on $n$ strands.
Let $U_k$ denote the set of indices of wires of $\W$ which 
lie above the vertex $v_k$ in $\pr_1^{-1}(y_k)$, and let $J_k
= (\overline V_k \setminus V_k) \cap U_k$.  Then the conjugating braids may 
be written as $\d _{k}=\prod_{j<i} A_{j,i}$, where the product is over all 
$i\in V_{k}$ and $j\in J_k$, see~\cite{CF}, \cite{CS3}.

Define a homomorphism 
$\boldsymbol{\Theta}_2(\mu):C_2\to C_2$ by
\begin{equation*}  \label{eq:Thetamu}
\boldsymbol{\Theta}_2(\mu) (e_{i}\wedge e_{j})
=\begin{cases}\Theta _{2}(\mu_{V_k})(e_{i}\wedge e_{j}) 
              &\text{ if } \{i,j\}\subset V_{k}\\
              e_{i}\wedge e_{j} &\text{ otherwise}.
\end{cases}
\end{equation*}
It is readily seen that $\boldsymbol{\Theta}_2(\mu)$ is invertible.  
Similarly, define $\boldsymbol{\Theta}_2(\d):C_2\to C_2$.  
A computation shows that 
$\d_k(t_i) = t_{J_k^i}\cdot t_i\cdot t_{J_k^i}^{-1}$, and that
$\Theta (\d_k)(e_{i})=
(1-t_{i})\cdot \nabla _{J_{k}^{i}} + t_{J_{k}^{i}}\cdot e_{i}$,  
for $i \in V_k$.    Thus, $\boldsymbol\Theta_2(\d)$ is also invertible.  
Proceeding as in~\ref{subsec:BAV}, we obtain the following.

\begin{thm}\label{thm:alexreal} 
The Alexander invariant of a complexified real arrangement
$\A $ has presentation
\begin{equation*}  \label{eq:alexreal}
C_3
\xrightarrow{\Delta'}  
L_0
\rightarrow B(\A )
\rightarrow 0,
\end{equation*}
where $L_0$ is the complementary summand to 
$K'_0=\oplus_V C_2'(V)$ in $C_2$, 
$\pi_0:C_2\to L_0$ is the canonical projection, 
and $\Delta' = \pi_0\circ \boldsymbol{\Theta}_2(\mu)^{-1} 
\circ \boldsymbol{\Theta}_2(\d)^{-1}\circ d_{3}$.
\end{thm}
Note that this presentation has only $\binom{n}{2}-b_2(M)$ 
generators, and $\binom {n}{3}$ relations.

\begin{rem} For an arbitrary complex arrangement, the map 
$\boldsymbol\Theta_2(\d)$ need not be invertible.  Thus the 
simplification of the presentation of the Alexander invariant 
afforded by the above result may not be available.
However, for any arrangement, we obtain an analogous 
simplified presentation for the $I$-adic completion, $\hB(\A)$, 
of the Alexander invariant of $\A$ in Corollary~\ref{cor:hBpres}.
\end{rem}

\section{Decomposition of the Alexander Invariant} 
\label{sec:DecompAlex}
We now relate the Alexander invariant of an arrangement $\A$ 
to a ``combinatorial'' Alexander invariant, determined by the
intersection lattice of $\A$.  For these purposes, we restrict our 
attention to central arrangements and their generic sections. 
It is enough to consider an affine arrangement, 
$\A=\{H_1,\dots,H_n\}$, of $n$ lines in $\C^2$ 
that is transverse to infinity. 
Recall that we identify set of rank two elements in the lattice of $\A$ 
and the collection of vertex sets of $\A$:  $L_2(\A)=\{V_1,\dots,V_s\}$.

\subsection{The Coarse Combinatorial Alexander Invariant} 
\label{subsec:ChainMap}
For each  
$V\in L_2(\A)$, let $G_V$ be the corresponding 
vertex group, and $B_V$ the corresponding Alexander invariant. 
Define the {\it coarse combinatorial Alexander invariant} of $\A$ to be 
\begin{equation*} \label{eq:localB}
\locB (\A)=\oplus_V B_V.
\end{equation*}
Notice that the module $B_V$ depends only on the cardinality $|V|$ of the
vertex set $V$.  Consequently, the module $\locB$ depends only on the 
number and
multiplicities of the elements of $L_2(\A)$.
This $\L$-module admits a free resolution 
\begin{equation*} \label{eq:localres}
\dots\rightarrow L_2 \xrightarrow{D_2} L_1
\xrightarrow{D_1} L_0 \rightarrow \locB \rightarrow 0
\end{equation*}
obtained by taking the direct sum of the resolutions \eqref{eq:BVres}: 
$L_k=\oplus_V C_2(V')\wedge C_{k}$, $D_k=\oplus_V\Delta_V^{k}$.  
(Since $\A$ is assumed to be transverse to infinity, 
$L_0=\oplus_V C_2(V')$ is indeed the complementary summand of 
$K'_0=\oplus_V C_2'(V)$ in $K_0 = C_2 = \oplus_V C_2(V)$.)

Let $\Psi_{V,\bullet}:C_{\bullet} \to C_2(V')\wedge C_{{\bullet}-2}$ 
be the chain map introduced in \eqref{eq:ChainV}.  
Define a chain map $\Psi_{\bullet}:C_{\bullet} \to L_{{\bullet}-2}$ 
by $\Psi_k = \sum_V\Psi_{V,k}$, for $k\ge 2$. 
\begin{prop} \label{prop:ups21}
The image of the composition 
$\Psi_2\circ \Phi: \bigoplus_V C_1(V') \to L_0$ 
is contained in the image of the map $D_1: L_1 \to L_0$. 
Therefore, there exists a map $\Gamma:\bigoplus_V C_1(V') \to L_1$ 
such that $D_1 \circ \Gamma = \Psi_2 \circ \Phi$.
\end{prop}
\begin{proof}
Let $A_V^\delta$ be a braid monodromy generator of $\A$, where 
$V=\{i_1,\dots,i_r\}$ and $\delta$ is some pure braid.  Using the pure braid
relations to rewrite $\delta$ if necessary, we may assume that this pure braid 
is a word in the generators $\{A_{r,s} \mid \{r,s\} \not\subset V\}$.

For $j\in V'$, we have $\Psi_2\circ \Phi(e_j) = 
\Psi_2 \circ \Theta_2(\delta)(\nabla_V \wedge e_j)$.  
Since $\delta\in\IA(F_n)$,
we have $\im(\Theta_k(\delta)-\id)\subset I\cdot C_k$.  
Hence, $\Theta_2(\d)(e_i\wedge e_j) = e_i\wedge e_j + W^\d_{i,j}$,
where $W^\d_{i,j}=\sum w_{p,q} e_p\wedge e_q$, with $w_{p,q}\in I$.
Thus $\Psi_2\circ \Phi(e_j) = 
\Psi_2(\nabla_V \wedge e_j) + \sum_{i\in V} t_{V^i}\cdot \Psi_2(W^\delta_{i,j})$.
Since $\Psi_2 (\nabla_V \wedge e_j)=\pi_{V'}\circ\Theta_2(\mu_V)^{-1}
(\Theta(\mu_V)(e_{i_1}\wedge e_j))=\pi_{V'}(e_{i_1}\wedge e_j)=0$, 
it suffices to show that
\begin{equation} \label{eq:imups1}
\Psi_2(W^\d_{i,j}) \in \im(D_1).
\end{equation}

For a vertex set $U \in L_2(\A)$, recall the natural projection
$\pi_U:C_k\to C_k(U)$, and denote by $I_U^\perp$ the ideal in
$\L$ generated by $\langle 1-t_k \mid k \notin U \rangle$.

\begin{claim} For each vertex set $U \in L_2(\A)$, we have
$\pi_U(W^\d_{i,j}) \in I_U^\perp \cdot C_2(U)$.
\end{claim}

Before proving this claim, let us show that it implies \eqref{eq:imups1}.
For $1\le p<q\le n$, let $V(p,q)$ 
denote the unique vertex set of $\A$ with $p,q\in V(p,q)$.  
If $w\cdot e_p\wedge e_q$ is a summand of $W^\d_{i,j}$, 
write $U=V(p,q)$.  Then, by the claim, we
have $w\in I_{U}^\perp$.  Now 
$\Psi_2(e_p\wedge e_q) = \Psi_{U,2}(e_p\wedge e_q)
\in C_2(U')$, and it is readily checked that
$I_{U}^\perp \cdot C_2(U') \subset \im(\Delta_U)$.
It follows that $\Psi_2(w\cdot e_p\wedge e_q) \in \im(\Delta_U)$.

Thus it suffices to prove the claim.
This may be accomplished by induction on the length of the word $\delta$.   
If $\delta = 1$, then $W^\d_{i,j}=0 $, and there is nothing to prove.   
If $\delta=A_{r,s}^{\pm 1}$, a computation shows that 
$\Theta_2(A_{r,s})(e_i\wedge e_j) = e_i\wedge e_j + W^{r,s}_{i,j}$ 
and $\Theta_2(A_{r,s}^{-1})(e_i\wedge e_j) = 
e_i\wedge e_j -(t_rt_s)^{-1} W^{r,s}_{i,j}$, where
\begin{equation*} \label{eq:T2rs}
W^{r,s}_{i,j}=
\begin{cases}
t_r(t_s-1)e_i\wedge e_r+t_r(1-t_r)e_i\wedge e_s 
& \text{if 
$i<j=r<s$,}\\
(1-t_j)[(1-t_s)e_i\wedge e_r+(t_r-1)e_i\wedge e_s] 
& \text{if 
$i<r<j<s$,}\\
(t_r-1)e_i\wedge e_s+(1-t_s)e_i\wedge e_r 
& \text{if 
$i<r<j=s$,}\\
(t_rt_s-1)e_r\wedge e_s 
& \text{if 
$i=r<j=s$,}\\
(t_r-1)e_i\wedge e_s+(t_s-1)[e_r\wedge e_i+(t_i-1)e_r\wedge e_s] 
& \text{if 
$r<i<j=s$,}\\
(1-t_i)[(1-t_s)e_r\wedge e_j+(t_r-1)e_s\wedge e_j] 
& \text{if 
$r<i<s<j$,}\\
t_r(t_s-1)e_r\wedge e_j+t_r(1-t_r)e_s\wedge e_j 
& \text{if 
$r=i<s<j$,}\\
(t_r-1)e_s\wedge e_j+(1-t_s)e_r\wedge e_j 
& \text{if 
$r<i=s<j$,}\\
t_r(t_s-1)e_r\wedge e_j+(t_r-1)[t_re_j\wedge e_s+(1-t_j)e_r\wedge e_s] 
& \text{if 
$r=i<j<s$,}\\
(t_j-1)[(1-t_s)e_r\wedge e_i+(1-t_r)e_i\wedge e_s] 
& \text{if 
$r<i<j<s$,}\\
\quad + (1-t_i)[(1-t_s)e_r\wedge e_j+(1-t_r)e_j\wedge e_s] \\
0 & \text{otherwise.}
\end{cases}
\end{equation*}
If $\{r,s\}\not\subset V$, it is readily checked that $\pi_U(W^{\d}_{i,j}) \in 
I_U^\perp\cdot C_2(U)$ for each $U \in L_2(\A)$.

In general, write $\d$ as the product of $A_{r,s}^{\pm 1}$ and $\d'$, and 
assume inductively that 
$\Theta_2(\d')(e_i\wedge e_j) = e_i\wedge e_j +W'_{i,j}$
satisfies $\pi_U(W'_{i,j}) \in I_U^\perp\cdot C_2(U)$ for each $U$. Then
$\Theta_2(\d)(e_i\wedge e_j)=\Theta_2(A_{r,s}^{\pm 1})(e_i\wedge e_j)+
\Theta_2(A_{r,s}^{\pm 1})(W'_{i,j})$, and by the above, it remains to
analyze the latter summand.  If $w\cdot e_p\wedge e_q$ is a summand of 
$W'_{i,j}$, then $w\in I_{V(p,q)}^\perp$ by induction.  
Case-by-case analysis then shows that each summand 
$x\cdot e_\ll\wedge e_m$
of $\Theta_2(A_{r,s}^{\pm 1})(w\cdot e_p\wedge e_q)$ satisfies
$x \in I_{V(\ll,m)}^\perp$.
This completes the proof of the claim, and hence that of the proposition.
\end{proof}

We can now formulate the main result of this subsection.

\begin{thm} \label{thm:BtoBlocal}
There exists a chain map $\Upsilon_\bullet$ from the presentation 
$K_{\bullet} \to B(\A)$ to the resolution $L_{\bullet}\to \locB (\A)$, 
\begin{equation*}\label{eq:chainmap}
\begin{CD}
&&&& K_1       @>\Delta>>        K_0@>>> B@>>>0\\
&&&& @VV{\Upsilon_1}V   @VV{\Upsilon_0}V @VV{\Pi}V \\
\dots@>>> L_2  @>D_2>> L_1  @>D_1>> L_0 @>>> \locB @>>> 0,
\end{CD}
\end{equation*}
given by $\Upsilon_0=\Psi_2$, and 
$\Upsilon_1(x,y) = \Gamma(x) + \Psi_3(y)$.  Furthermore, 
the resulting map $\Pi:B\to \locB$ is surjective.  
\end{thm}
\begin{proof}  It is immediate from the above that $\Upsilon_\bullet$ is
a chain map.  Thus it suffices to show that the map 
$\Upsilon_0:K_0\to L_0$, which by definition equals
$\Psi_2=\sum_V\pi_{V'}\circ\Theta_2(\mu_V)^{-1}: 
C_2\to \oplus_V C_2(V')$,
is surjective.
With respect to the decomposition $K_0=K_0'\oplus L_0
=\left( \oplus_V C_2'(V)\right) \oplus 
\left( \oplus_V C_2(V')\right)$, we have 
$\Upsilon_0=\pmatrix \Upsilon'_0 & \id_{L_0}\endpmatrix^{\top}: 
K_0'\oplus L_0\to L_0$, where   
$\Upsilon'_0(e_i\wedge x) = (e_i-\nabla_V)\wedge x$, 
for $i=\min V$ and $x \in C_1(V')$.  
Thus $\Upsilon_0:K_0\to L_0$ is surjective.  
\end{proof}

\subsection{Decomposition of the Completion}
\label{subsec:Decomp}
Recall that if $B$ is a $\L$-module, then $\hB$ denotes its 
$I$-adic completion, and that if $f:A\to B$ is a map of $\L$-modules, 
we write $\widehat f:\widehat A \to \hB$ for the extension of $f$ to 
the completions.  The $I$-adic completion functor takes chain complexes 
to chain complexes, and chain maps to chain maps.  

\begin{thm} \label{thm:decomp}
The chain map 
$\widehat\Upsilon_\bullet:\widehat K_\bullet\to \widehat L_\bullet$ 
induces an isomorphism $\hB \xrightarrow{\sim}\hlocB$ if and only if 
the map $\widehat\Psi_3:\widehat C_3 \to \widehat L_1$ is surjective.
\end{thm}
\begin{proof}  
Consider the mapping cone, $K_\bullet(\widehat \Upsilon)$, of 
$\widehat\Upsilon_\bullet$, given by
\begin{equation} \label{eq:mapcone}
\dots \to \widehat L_2 \oplus \widehat K_1 \xrightarrow{\partial_2}
\widehat L_1 \oplus \widehat K_0 \xrightarrow{\partial_1} \widehat L_0,
\end{equation}
where $\partial_2(x,y) = (\widehat{D}_2(x) - 
\widehat \Upsilon_1(y),\widehat\Delta(y))$ and 
$\partial_1(x,y) = \widehat{D}_1(x) + \widehat\Upsilon_0(y)$, 
and the short exact sequence of chain complexes
\begin{equation} \label{eq:seq}
0 \to \widehat L_\bullet \xrightarrow{\iota_\bullet} 
K_\bullet(\widehat \Upsilon) \xrightarrow{\pi_\bullet} 
\widehat K_{\bullet - 1} \to 0,
\end{equation}
where $\iota_\bullet$ and $\pi_\bullet$ denote the natural inclusion 
and projection.  Since $\widehat\Upsilon_0$ is surjective, we have 
$H_0(K_\bullet(\widehat \Upsilon))=0$.  Also, since $\widehat L_\bullet$ 
is a resolution, $\widetilde H_*(\widehat L_\bullet)=0$.  
Thus the associated long exact sequence in homology reduces to
\begin{equation} \label{eq:Hseq}
0 \to H_1(K_\bullet(\widehat \Upsilon)) \to H_0(\widehat K_\bullet)
\xrightarrow{\widehat\Upsilon_{*}} H_0(\widehat L_\bullet) \to 0.
\end{equation}
The map $\widehat\Upsilon_{*}:H_0(\widehat K_\bullet) = 
H_0(\widehat L_\bullet)$ identifies canonically with 
$\widehat\Pi:\hB \to\hlocB$.  
Thus, it suffices to show that $\coker(\widehat\Psi_3)=0$ 
if and only if $H_1(K_\bullet(\widehat \Upsilon))=0$.

Recall the map $\Phi:\oplus_V C_1(V') \to C_2$ from 
Theorem~\ref{thm:alexarr}.  
Recall also (from the proof of Proposition~\ref{prop:ups21}) 
that $\Phi|_{C_1(V')}=\Phi_V\mod I$.  Thus, 
$\widehat\Phi|_{\widehat C_1(V')} = \widehat\Phi_V \mod \mathfrak{m}$.  
Using the identification $C_1(V')\xrightarrow{\sim} C'_2(V)$, 
$x\mapsto e_i\wedge x$, where $i=\min V$, and the 
projection onto first factor $p':K_0=K_0'\oplus L_0\to K_0'$, define 
$\Phi':= p' \circ \Phi : K'_0 \to  K'_0$.  
Since $\widehat\nabla_V=\sum_{i\in V} e_i$, the map 
$\widehat p'\circ \widehat\Phi_V:\widehat C_1(V')\to \widehat C'_2(V)$ 
coincides with the above identification.   
Hence, $\widehat\Phi'=\id\mod \mathfrak{m}$.  Consequently, 
$\widehat\Phi'$ is an isomorphism.

We now alter the short exact sequence~\eqref{eq:seq}.  Write
$K_1(\widehat\Upsilon) = \widehat L_1 \oplus 
\widehat K'_0 \oplus \widehat L_0$ 
and $K_2(\widehat\Upsilon) = \widehat L_2 \oplus 
\widehat K'_0 \oplus \widehat C_3$, 
and define $\rho \in \Aut K_1(\widehat\Upsilon)$ and
$\psi \in \Aut K_2(\widehat\Upsilon)$ by
\begin{align*}
\rho(x,y,z) & = \left( x + \widehat\Gamma \circ
\widehat{\Phi'}^{-1}(y), y, D_1(x)+\widehat\Upsilon'_0(y) +z\right),\\
\psi(x,y,z) &= \left( x, \widehat{\Phi'}^{-1}(y) - \widehat{\Phi'}^{-1} 
\circ p'\circ\widehat d_3(z), z\right). 
\end{align*}
Note that the restriction of $\psi$ to $\widehat L_2$
is the identity, and let $\bar\psi$ denote the restriction of $\psi$ to
$\widehat K_1 = \widehat K'_0 \oplus \widehat C_3$.  
We modify the sequence~\eqref{eq:seq}
as indicated below.
\begin{equation}\label{eq:MapConeSeq}
\begin{CD}
\dots @>>> \widehat L_2 @>\widehat D_2>>  \widehat L_1
@>\widehat D_1>>  \widehat L_0\\
&&@VV{\psi^{-1} \circ \iota_2}V @VV{\rho \circ \iota_1}V   
@VV{\iota_0=\id}V \\
\dots @>>> \widehat L_2 \oplus \widehat K_1 
@>{\rho \circ \partial_2 \circ \psi}>>  
\widehat L_1 \oplus \widehat K_0
@>{\partial_1 \circ \rho^{-1}}>>  \widehat L_0\\
&&@VV{\pi_2}V @VV{\pi_1 \circ \rho^{-1}}V  @VVV\\
&&\widehat K_1 
@>{\widehat\Delta \circ \bar\psi}>>  
\widehat K_0 @>>> 0
\end{CD}
\end{equation}
Since $\bar\psi\circ\pi_2 = \pi_2\circ\psi$, this diagram commutes.

Consider the map $\Xi:\widehat L_2 \oplus \widehat C_3
\to \widehat L_1$ defined by
\begin{equation} \label{eq:H1pres}
\Xi(x,z)=\widehat D_2(x)+ \widehat\Gamma \circ
\widehat{\Phi'}^{-1} \circ p' \circ \widehat d_3(z) - 
\widehat\Psi_3(z).
\end{equation}
Computations with the definitions (making use of the fact that 
$\widehat\Upsilon_\bullet$ is a chain map) reveal that
$\partial_1 \circ \rho^{-1}(x,y,z)=z$ and 
$\rho \circ \partial_2 \circ \psi(x,y,z)=(\Xi(x,z),y,0)$.  
Thus, $\Xi$ provides a presentation for the module 
$H_1(K_\bullet(\widehat \Upsilon))$, 
and this module is trivial if and only if $\Xi$ is surjective.
Since $\im(\widehat D_2)\subset  {\mathfrak m} \cdot \widehat L_1$
and $\im(\widehat\Gamma \circ
\widehat{\Phi'}^{-1} \circ p' \circ \widehat d_3)\subset  
{\mathfrak m} \cdot \widehat L_1$, the map $\Xi$ is surjective 
if and only if $\widehat\Psi_3$ is surjective.
\end{proof}

The above proof has several consequences, even in the instance when the
map $\widehat\Psi_3$ is not surjective, see below and
Theorem~\ref{thm:theta3}.   These results hold for an arbitrary
arrangement (real or complex,  compare~\ref{thm:alexarr}
and~\ref{thm:alexreal}) that is transverse  to infinity.

\begin{cor} \label{cor:hBpres}
The $I$-adic completion $\hB$ of the Alexander invariant
an arrangement $\A$ has a presentation
with $\binom {n}{2}-b_2(M(\A))$ generators, and 
$\binom {n}{3}$ relations.
\end{cor}
\begin{proof}
From the commutative diagram~\eqref{eq:MapConeSeq}, we have the
presentation $\widehat\Delta \circ \bar\psi:\widehat K_1 \to  \widehat
K_0$ for $\hB$.   Let $p'':K_0 = K'_0 \oplus  L_0 \to L_0$ denote the
projection  onto second factor.   With respect to the decompositions 
$\widehat K_1 = \widehat K'_0 \oplus \widehat C_3$ and 
$\widehat K_0 = \widehat K'_0 \oplus \widehat L_0$, the map
$\widehat\Delta \circ \bar\psi:\widehat K_1\to \widehat K_0$ is given by
\begin{equation*} \label{eq:deltapsi}
\widehat\Delta \circ \bar\psi(x,y) = \left(x,
\widehat p''\left(\widehat\Phi\circ\widehat{\Phi'}^{-1}(x) - 
\widehat\Phi\circ\widehat{\Phi'}^{-1}\circ\widehat p'\circ \widehat d_3(y) +
\widehat d_3(y) \right)\right). 
\end{equation*}
Define 
\begin{equation} \label{eq:deltasharp}
\widehat\Delta^\sharp=
\widehat p''\circ \left( \id - \widehat \Phi\circ \widehat {\Phi'}^{-1} 
\circ \widehat p'\right)\circ \widehat d_3:\widehat C_3 \to \widehat L_0,
\end{equation}
and define $\chi \in \Aut \widehat K_0$ by 
$\chi(x,y) = (x,y - \widehat p''\circ\widehat\Phi\circ \widehat{\Phi'}^{-1}(x))$.  
Then $\chi\circ \widehat\Delta \circ \bar\psi(x,y) = (x, \widehat\Delta^\sharp(y))$.  
Thus, $\widehat\Delta^\sharp$ provides a presentation for $\hB$ with 
the specified numbers of generators and relations.
\end{proof}

\section{Combinatorics and the Chen Groups}
\label{sec:Combinatorics}
In this section, we examine the relationship between the results obtained in
the previous sections and the combinatorics of the arrangement $\A$.  

An invariant of $\A$  is called {\it combinatorial} if it is determined by the 
isomorphism type of the lattice $L(\A)$.  As is well-known from \cite{F1}, 
the ranks $\phi_k$ of the LCS quotients of the group of $\A$ are combinatorially 
determined.  Thus,  the ranks of the first three Chen groups of $\A$ 
are combinatorial.  We now describe some explicit combinatorial bounds 
and formulas for the ranks $\theta_k$ of the Chen groups of $\A$.    

\subsection{A Bound on Chen Ranks} 
\label{subsec:ChenBound}
Recall that the coarse combinatorial Alexander invariant $\locB$ 
of $\A$ is the direct sum $\oplus_V B_V$ of the Alexander invariants 
of the vertex groups $G_V$, indexed by $V\in L_2(\A)$, the rank two 
elements of $L(\A)$.  For $k\ge 2$, define the 
{\it coarse combinatorial Chen ranks} by 
\begin{equation*} \label{eq:thetacc}
\theta^\loc_k(\A) = \sum_{V\in L_2(\A)} \theta_k(G_V).
\end{equation*}
From Theorem~\ref{thm:massey} and Theorem~\ref{thm:BtoBlocal}, 
we obtain the following.
\begin{cor}\label{cor:chenlowerbound}
For $k\ge 2$, the ranks of the Chen groups of $\A$ are bounded below
by the coarse combinatorial Chen ranks:  $\theta_k(\A) \ge \theta^\loc_k(\A)$.
\end{cor}

To compute these lower bounds explicitly, recall that 
$G_V\cong F_{r-1} \times \Z^{n-r}$, where $r=|V|$.  By 
Corollary~\ref{cor:chenprod}, the ranks of the Chen groups 
of $G_V$ are given by 
$\theta_k(G_V) = (k-1) \binom{k+r-3}{k}$ for $k \ge 2$.  
Let $c_r$ denote the number of elements of $L_2(\A)$ of multiplicity $r$,
and write $\binom{m}{k}=0$ if $m<k$.  Then,
\begin{equation*} \label{eq:LocalComb}
\theta_k^\loc =\sum_{V \in L_2(\A)}  (k-1) \binom{k+|V|-3}{k}
= \sum_{r\ge 3} c_r (k-1) \binom{k+r-3}{k}, 
\end{equation*}
and so $\theta_k^\loc$ is determined by (only) the multiplicities of the elements 
of $L_2(\A)$.

\begin{rem} The ranks of the lower central series quotients of the group
$G$ of $\A$ satisfy analogous lower bounds: $\phi_k(\A) \ge \phi^\loc_k(\A)
=\sum_{V \in L_2(\A)} \phi_k(G_V)$, see~\cite{F2}~Proposition~3.8.
\end{rem}

\begin{rem} \label{rem:mobius}
The lower bounds for the ranks of the Chen groups of $\A$ 
may be expressed in terms of the M\"obius function $\mu:L(\A)\to\Z$: 
$\theta_k(\A) \ge \sum_{V\in L_2(\A)} (k-1)\binom{k+\mu(V)-2}{k}$.
\end{rem}

We now analyze the difference $\theta_k-\theta_k^\loc$.  
Recall that 
$\theta_k(\A) = \rank({\mathfrak m}^{k-2}\hB/{\mathfrak m}^{k-1}\hB)$.
Checking that the image of the map 
$\widehat\Delta^\sharp:\widehat C_3 \to \widehat L_0$  
defined in~\eqref{eq:deltasharp} is contained
in ${\mathfrak m} \cdot \widehat L_0$, we see that 
$\theta_2(\A) = \rank(\hB/{\mathfrak m}
\hB) = \binom {n}{2}-b_2(M(\A))$, 
and thus $\theta_2(\A)=\theta^\loc_2(\A)$.  

Recall the mapping cone $K_{\bullet}(\widehat\Upsilon)$ 
from~\eqref{eq:mapcone}, and set 
$H=H_1(K_\bullet(\widehat \Upsilon))$.  

\begin{thm} \label{thm:theta3}
For $k \ge 3$, the rank of the $k^{\text{th}}$ Chen group of $\A$ 
is given by
$\theta_k(\A) = \rank({\mathfrak m}^{k-3}H/{\mathfrak m}^{k-2}H)
+ \theta^\loc_k(\A)$.  In particular, 
$\theta_3(\A) = \rank (\coker\widehat\Psi_3)  + \theta^\loc_{3}(\A)$. 
\end{thm}
\begin{proof}
Consider the short exact sequence~\eqref{eq:Hseq}, rewritten as
\begin{equation}  \label{eq:PiB}
0 \to H 
\xrightarrow{\tau} \hB
\xrightarrow{\widehat\Pi} \hlocB \to 0.
\end{equation}
Altering the commutative diagram~\eqref{eq:MapConeSeq} 
using the isomorphism $\chi \in \Aut \widehat K_0$ defined 
in the proof of Corollary~\ref{cor:hBpres}, 
we see that $\tau$ is induced by 
$\chi \circ \pi_1 \circ \rho^{-1}:
\widehat L_1 \oplus \widehat  K'_0 \oplus \widehat L_0 \to 
\widehat K_0$, $(x,y,z) \mapsto (y, z - \widehat D_1(x))$.  
Thus the restriction of this map to
$\ker(\partial_1 \circ\rho^{-1})=  \widehat L_1 \oplus \widehat K'_0$ 
is given by $(x,y) \mapsto (y, - \widehat D_1(x))$.
Since $\hB = \coker(\chi\circ \widehat\Delta \circ \bar\psi)$ and
$\chi\circ \widehat\Delta \circ \bar\psi(y,0) = (y,0)$,
the map $\tau:H\to \hB$ in homology is induced by 
$\widehat D_1:\widehat L_1 \to \widehat L_0$.
Since  
$\im(\widehat D_1(z)) \subset{\mathfrak m} \cdot \widehat L_0$, 
from the exact sequence~\eqref{eq:PiB} we have
\begin{equation*}
\rank({\mathfrak m}^{\ll}\hB/{\mathfrak m}^{\ll-1}\hB) = 
\rank({\mathfrak m}^{\ll-1}H/{\mathfrak m}^{\ll-2}H)
+ \rank({\mathfrak m}^{\ll} \hlocB/{\mathfrak m}^{\ll-1}\hlocB).
\end{equation*}
It follows from Theorem~\ref{thm:massey} 
that the ranks of the Chen groups of $\A$ are as asserted.

In particular, the third Chen group of $\A$ has rank
$\theta_3(\A) = \rank(H/{\mathfrak m} H) + \theta^\loc_3(\A)$.
Recall the presentation,
$\Xi:\widehat L_2 \oplus \widehat C_3 \to \widehat L_1$,
for the module $H$ from~\eqref{eq:H1pres}.  
Using elementary row and column operations, 
we obtain a presentation $\Xi':\hL^a \to \hL^b$ from this with 
$b=\rank \widehat L_1 - \rank\widehat\Psi_3$ generators.
Checking that $\im(\Xi') \subset {\mathfrak m}\cdot \hL^b$, 
we get 
$\rank(H/{\mathfrak m} H) = \rank(\coker \widehat\Psi_3)$.
\end{proof}

\subsection{Decomposition is Combinatorial} 
\label{subsec:theta3}
Let $\hat\epsilon:\hL \to \Z$ be the augmentation map, which takes a power 
series to its constant coefficient.  If $\widehat F=\hL^p$ is a free module, 
denote its image under $\hat\epsilon$ by $\overline F=\Z^p$, and if 
$\widehat f:\widehat F\to\widehat F'$ is a $\hL$-linear map, denote 
its image by $\overline f:\overline F \to \overline F'$.   

\begin{lem}
The rank of $\widehat f$ is equal to the rank of $\overline f$.
\end{lem}

\begin{proof}
Suppose the rank of $\overline f:\Z^p \to \Z^q$ is $r$.  Then there are 
integral matrices $X \in \GL(p,\Z)$ and $Y \in \GL(q,\Z)$ so that 
$X \cdot \overline f \cdot Y = \begin{pmatrix} I_r&0\\0&0 \end{pmatrix}$, 
where $I_r$ denotes the $r\times r$ identity matrix.  By definition, 
$\widehat f = \overline f + Z$, where the entries of $Z$ are in $\mathfrak m$.
Thus, $X \cdot \widehat f \cdot Y = X \cdot \overline f \cdot Y+
X \cdot Z \cdot Y$.  Clearly, the rank of $X \cdot \widehat f \cdot Y$, 
and hence that of $\widehat f$, is $r$.

The converse follows from the functoriality of the construction. 
\end{proof}

We now show that the rank of the map 
$\widehat\Psi_3:\widehat C_3\to \widehat L_1$ is 
combinatorially determined.  Thus, the criterion for decomposition 
of the $I$-adic completion of the Alexander invariant of 
Theorem~\ref{thm:decomp}---the surjectivity of $\widehat\Psi_3$---is 
combinatorial as well.  By the lemma, it suffices to show that the 
rank of $\overline \Psi_3:\overline C_3\to \overline L_1$ is 
combinatorially determined. 

For this, let $\A$ and $\A^*$ be lattice-isomorphic arrangements 
of $n$ lines in $\C^2$ (which are transverse to the line at infinity). 
Let $\W$ be a braided wiring diagram associated to
$\A$, and let $\overline \Psi_3:\overline C_3\to \overline L_1=
\bigoplus_V \overline C_2(V') \wedge \overline C_1$ 
be the map defined by the vertex sets $\{V_1,\dots,V_s\}$
of $\W$.  Choose arbitrary orderings of the hyperplanes and rank two lattice
elements of $\A^*$, and denote the elements of $L_2(\A^*)$ by
$\{U_1,\dots,U_s\}$.  Then formally construct the map 
$\overline\Psi^*_{3} = \sum_U \overline\Psi_{U,3}:
\overline C_3 \to \overline L_1^*=
\bigoplus_U \overline C_2(U') \wedge \overline C_1$ 
using~\eqref{eq:ChainV}, 
the Magnus embedding, and the augmentation map~$\hat\epsilon$.

Since $\A$ and $\A^*$ are lattice-isomorphic, there are permutations 
$\omega\in\Sigma_n$ and $\nu\in\Sigma_s$ so that $\omega(V_k) =
U_{\nu(k)}$ for each $k$, $1\le k\le s$.  The permutation $\omega$ 
induces an isomorphism $\omega_k: C_k\to C_k$ defined by 
$\omega_k(e_J) = e_{\omega(J)}$.  The map $\overline\Psi_3$ is
combinatorially determined in the sense of the following.

\begin{prop} There is an isomorphism
$\overline\xi : \overline L_1 \to \overline L_1^*$ so that 
$ \overline\xi \circ \overline{\Psi}_{3} = 
\overline\Psi^*_{3} \circ \overline\omega_3$.
\end{prop}
\begin{proof}  Let $V$ be a vertex set of $\A$, and 
$U = \omega(V)$ be the corresponding vertex set of $\A^*$.
Define a map $\xi^U_V:C_2(V') \wedge C_1 \to C_2(U') \wedge C_1$ by
\begin{equation} \label{eq:latticeiso}
\xi^U_V=(\pi_{U'}\wedge\id) \circ \Theta_2(\mu_U)^{-1}\circ \omega_3 
\circ \Theta_2(\mu_V)\circ (\iota_{V'}\wedge \id).  
\end{equation}
Clearly,  $\xi^U_V$ is an isomorphism, with inverse $\xi^V_U$.  
Moreover, $\xi^U_V \circ\Psi_{V,3} = 
\Psi_{U,3} \circ \omega_3$.    The collection $\{\xi^V_U\}$ defines a 
map $\xi :  L_1 \to L_1^*$, which yields the desired isomorphism 
$\overline\xi$.
\end{proof}

Combining these results with those of the previous section, we obtain

\begin{thm} \label{thm:ChenComb}
The rank of the third Chen group of the arrangement $\A$ is 
given by the combinatorial formula
$\theta_3(\A) = \rank(\coker\overline\Psi_3) + \theta^\loc_3(\A)$.

Furthermore, if the map $\overline \Psi_3:\overline C_3\to \overline L_1$  
is surjective, then the $I$-adic completion of the Alexander invariant of $\A$ 
decomposes as a direct sum: 
$\hB \cong \hlocB=\oplus_V \hB_V$, and the ranks of the Chen 
groups of $\A$ are given by $\theta_k(\A)=\theta_k^\loc(\A)$ for all $k\ge 2$.
\end{thm}

\begin{rem} If $\A$ is an arrangement for which the map
$\overline \Psi_3:\overline C_3\to \overline L_1$ is {\em not} surjective,
the ranks $\theta_k(\A)$ of the Chen groups of $\A$ for $k \ge 4$ may 
be computed 
using the Groebner basis algorithm described in~\ref{subsec:Chen}.
Alternatively, in light of Theorem~\ref{thm:theta3}, one can apply this 
algorithm to the presentation~\eqref{eq:H1pres} of the module $H$
(or the smaller presentation described in the proof of~\ref{thm:theta3})
to determine $\rank({\mathfrak m}^{k-3}H/{\mathfrak m}^{k-2}H)=
\theta_k(\A) - \theta^\loc_k(\A)$.
\end{rem}

\section{Examples}
\label{sec:Examples}

In this section, we illustrate the results of the previous sections 
by means of several  explicit examples.  We
order the hyperplanes of an arrangement $\A=\{H_1,\dots,H_n\}$ in the
order indicated by the defining polynomial $Q(\A) = \prod_{k=1}^n
\ll_k$ (so $H_k = \ker\ll_k$).

\begin{exm}\label{exm:6lines} Consider the central 3-arrangement $\A$
with defining polynomial 
$Q=xyz(y+z)(x-z)(2x+y)$.  Randell \cite{Ra} noted that this
arrangement is not $K(\pi ,1)$, 
and that there is no aspherical arrangement with
the same lattice in ranks one  and two.  Arvola \cite{Ar} further
showed that the group of this arrangement is not of  type FL.

The rank two elements of the lattice of $\A$ are
\begin{equation*}
L_2(\A) = \{\{1, 2, 6\}, \{1, 3, 5\}, \{2, 3, 4\}, \{1, 4\}, \{2,
5\}, \{4, 5\}, \{3, 6\}, \{4, 6\}, \{5, 6\}\}.
\end{equation*}
It is readily checked that the map $\overline\Psi_3:\Z^{20} 
\to \Z^{12}$ is surjective.  By Theorem~\ref{thm:ChenComb}, 
$\hB \cong\hlocB=\hB_{\{1,2,6\}} \oplus \hB_{\{1,3,5\}} 
\oplus \hB_{\{2,3,4\}}$.  It
follows that the ranks of the Chen groups of $\A $ are 
$\theta_{1}=6$ and $\theta _{k}=3(k-1)$ for $k\ge 2$.  
Notice that these ranks coincide with those of the Chen 
groups of a direct product of three free groups on two generators, 
though clearly $G \not\cong F_{2}\times F_{2}\times F_{2}$.

Using Theorem~\ref{thm:alexreal} and elementary row operations, 
one can show that the Alexander invariant itself decomposes as a
direct sum, $B \cong \locB=
B_{\{1,2,6\}} \oplus B_{\{1,3,5\}} \oplus B_{\{2,3,4\}}$.
\end{exm}

\begin{exm}\label{exm:braid} The braid arrangement $\A_4$ is the
smallest arrangement  
for which the completion
of the Alexander invariant does not decompose.   The polynomial
$Q=xyz(x-y)(x-z)(y-z)$ defines a central 3-arrangement whose
complement is homotopy equivalent to that of $\A_4$.  The rank two
elements of $L(\A_4)$ (the partition lattice) are
\begin{equation*}
L_2(\A_4) =\{\{1, 2, 4\}, \{1, 3, 5\}, \{2, 3, 6\}, \{3, 4\}, \{2, 5\},
\{4, 5, 6\}, \{1, 6\}\}.
\end{equation*}
The map $\overline\Psi_3:\Z^{20}\to \Z^{16}$
is {\it not} surjective.  Thus $\hB(\A_4)$ does not decompose.  
A basis for $\coker\overline\Psi_3$ is given by the two elements
\begin{align*}
\kappa_1&=e_{\{2,4\}}\wedge(e_6-e_3)+e_{\{3,5\}}\wedge(e_4-e_6)+
e_{\{3,6\}}\wedge(e_1-e_4)+e_{\{5,6\}}\wedge(e_3-e_1),\\
\kappa_2&=e_{\{2,4\}}\wedge(e_6-e_5)+e_{\{3,5\}}\wedge(e_2-e_6)+
e_{\{3,6\}}\wedge(e_1-e_5)+e_{\{5,6\}}\wedge(e_2-e_1).
\end{align*}
Since $\rank\overline\Psi_3=14$, we have $\theta_3(\A_4)=10$.
\end{exm}

\begin{rem}
Note that the rank of the third Chen group of $\A_4$ is equal to that of
the product arrangement defined by $xy(y-x)z(z-x)(z-2x)$.  In general,
by the LCS formula~\cite{FR1}, the ranks of the lower central series quotients 
of the pure braid group $P_n$ are equal to those of the direct product 
of free groups $\Pi_n = F_{n-1}\times\dots\times F_1$.  
These groups are distinguished by their Chen groups.  
For $k \ge 4$, we have $\theta_k(\Pi_n) =
(k-1)\binom{n+k-2}{k+1}$, by Example~\ref{exm:prodFn}, 
and $\theta_k(P_n) = (k-1)\binom{n+1}{4}$, 
by the main result of~\cite{CS1}.  
Thus, $\theta_k(P_n) \neq \theta_k(\Pi_n)$ for $n\ge 4$, 
and the groups $P_n$ and $\Pi_n$ are not isomorphic.
\end{rem}

\begin{rem}\label{rem:suba4} Example~\ref{exm:braid} provides an easy means 
for detecting when the completion of the Alexander invariant of an 
arrangement $\A$ does not  decompose.  If $\calS \subset \A$ is a 
subarrangement which is lattice-isomorphic to the braid arrangement 
$\A_4$, one can use the above elements of
$\coker\overline\Psi_3(\A_4)$ and maps of the
form~\eqref{eq:latticeiso} to generate non-trivial elements of 
$\coker \overline\Psi_3(\A)$.

It is interesting to note that (the matroid of) such an arrangement $\A$ 
has ``non-local decomposable relations,'' see~\cite{F3}.
\end{rem}

\begin{exm}[Diamond]\label{exm:diamond} Let $\mathcal D$ be the central
$3$-arrangement  with defining polynomial
$Q(\mathcal D)=x(x+y+z)(x+y-z)y(x-y-z)(x-y+z)z$.  
This is a free,  simplicial
arrangement for which the LCS formula does not hold, and the
Orlik-Solomon algebra is not quadratic, see \cite{F1},
\cite{FR2}, \cite{SY}.   The rank two elements of $L(\mathcal D)$ are
\begin{equation*}
L_2(\mathcal D) =\{\{3,4,5\}, \{1,2,5\}, \{1,4\}, \{1,3,6\}, 
\{2,4,6\}, \{1,7\}, \{2,3,7\}, \{4,7\}, \{5,6,7\}\}.
\end{equation*}

This arrangement has three distinct subarrangements
lattice-isomorphic to $\A_4$.  
One such subarrangement is
$\calS=\mathcal D \setminus \{H_1\}$.
Define $\omega:[6]\hookrightarrow [7]$ by
$1\mapsto 2$, $2\mapsto 3$, $3\mapsto 4$, $4\mapsto 7$, $5\mapsto 6$, 
$6\mapsto 5$.
The map $\omega$ gives rise to a lattice-isomorphism $L(\A_4)
\xrightarrow{\sim} L(\calS)$.  Define 
$\overline\xi:\bigoplus_{V\in L_2(\A_4)} \overline C_2(V')\wedge \overline C_1 
\to \bigoplus_{U\in L_2(\mathcal D)}  \overline C_2(U')\wedge \overline C_1$  
as in~\eqref{eq:latticeiso}:
$\overline\xi(e_{\{5,6\}}\wedge e_j) = -e_{\{6,7\}}\wedge 
\overline\omega_1(e_j)$ and $\overline\xi(e_K\wedge e_j) = 
\overline\omega_3(e_K \wedge e_j)$ 
for $K \neq \{5,6\}$. Then 
\begin{align*}
\overline\xi(\kappa_1)&=e_{\{3,7\}}\wedge(e_5-e_4)+
e_{\{4,6\}}\wedge(e_7-e_5)+
e_{\{4,5\}}\wedge(e_2-e_7)+e_{\{6,7\}}\wedge(e_2-e_4),\\
\overline\xi(\kappa_2)&=e_{\{3,7\}}\wedge(e_5-e_6)+
e_{\{4,6\}}\wedge(e_3-e_5)+
e_{\{4,5\}}\wedge(e_2-e_6)+e_{\{6,7\}}\wedge(e_2-e_3)
\end{align*}
are in the cokernel of $\overline\Psi_3: \Z^{35}\to\Z^{30}$.  
We obtain 6 distinct elements of~$\coker \overline\Psi_3$ in this way.  
However, there is a relation among them.
We have $\rank\overline\Psi_3=25$, and~$\theta_3(\mathcal D)=17$.

The ranks of the higher Chen groups 
may be found via the 
Groebner basis algorithm of~Theorem~\ref{thm:grB}.  
By Example~\ref{exm:cone}, we can 
simplify the computation 
by working with the decone of $\mathcal D$ defined by 
$Q(\mathcal D)|_{z=1}$.  Rotating this arrangement counterclockwise
to insure that first coordinate projection is generic, we obtain a
wiring diagram
$\W=\{\{3,4,5\},\{1,2,5\},\{1,4\},\{1,3,6\},\{2,4,6\}\}$.  
The image of the associated braid 
monodromy $\a :F_{5}\to P_{6}$ is 
generated by 
$\{A_{3,4,5}, A_{1,2,5}, A_{1,4}^{A_{3,4}}, 
A_{1,3,6}, A_{2,4,6}^{A_{3,4}A_{3,6}}\}$.  
From the presentation $\L ^{20} \xrightarrow{\Delta} \L ^{6}
\rightarrow B \rightarrow 0$  
provided by Theorem~\ref{thm:alexreal}, we find 
$\theta_k(\A)=9(k-1)$, for $k\ge 4$.  
\end{exm}

We have found a number of other arrangements for which the completion
of the Alexander invariant does not decompose as a direct sum.  For
example, for the Coxeter arrangement of type $\text{B}_3$, we have
$\overline C_3=\Z^{84}$, $\overline L_1=\Z^{85}$, and 
$\overline\Psi_3:\overline C_3\to \overline L_1$
is obviously not surjective.    
(This arrangement also has (many) subarrangements lattice-isomorphic 
to $\A_4$.)  More subtle examples include the following.

\begin{exm}[MacLane]\label{exm:maclane}  The polynomials
\begin{equation*}
Q^{\pm}=x y (y-x) z (z-x-\omega^2 y) (z+\omega y) (z-x) (z+\omega^2 x+\omega y)
\end{equation*}
where $\omega=(-1\pm \sqrt{-3})/2$, define complex conjugate
realizations $\A^\pm$ of the MacLane matroid (the $8_3$ configuration).  
These arrangements were used by Rybnikov \cite{Ry} in his 
construction of lattice-isomorphic arrangements with distinct 
fundamental groups.  Rybnikov's arrangements are not distinguished 
by their Chen groups.

Complex conjugation induces an isomorphism of the 
groups of $\A^+$ and $\A^-$,
and thus an isomorphism of the Alexander 
invariants,  $B^+\cong B^-$.   
Neither of these arrangements has a subarrangement
lattice-isomorphic to the braid arrangement $\A_4$.  
Nevertheless, the
maps $\overline\Psi_3^{\pm}$ are not surjective, and the modules $\hB^{\pm}$ 
do not decompose.   

The hyperplane $H=\{7x-5y+z=1\}$ is generic with 
respect to both these arrangements.  Moreover, 
the projection $\pr :\C ^{2}\to \C$ defined by $\pr (x,y)=13x-4y$ 
is generic with respect to both affine $2$-arrangements $H\cap\A^\pm$.
Changing coordinates accordingly, we obtain braided wiring diagrams
$\W^\pm=\{ V_1,\b_1^\pm,V_2,\b_2^\pm,\dots,\b_{11}^\pm,V_{12}\}$ 
with vertex sets and intermediary braids given by:
\begin{equation*}  \label{eq:bwire}
\matrix
V_{1}=\{3,4\},\hfill&\b_{1}^{\pm} = 1,\hfill&\qquad
V_{7\,\,\,}=\{2,5,7\},\hfill&\b_{7}^{\pm}=\s_3^{\mp 1},\hfill\\
V_2=\{3,5,6\},\hfill&\b_{2}^{\pm} = 1,\hfill&\qquad
V_{8\,\,\,}=\{6,7\},\hfill&\b_{8}^{\pm} = \s_4^{\pm 1}\s_3^{\pm 1},\hfill\\
V_3=\{3,7,8\},\hfill&\b_{3}^{\pm}=1,\hfill&\qquad
V_{9\,\,\,}=\{1,6,8\},\hfill&\b_{9}^{\pm}=\s_1^{\mp 1},\hfill\\
V_4=\{2,4,6\},\hfill&\b_{4}^{\pm}=\s_4^{\mp 1}\s_3^{\mp 1},\hfill&\qquad
V_{10}=\{1,4,7\},\hfill&\b_{10}^{\pm}=\s_2^{\mp 1},\hfill\\
V_5=\{2,8\},\hfill&\b_{5}^{\pm} = 1,\hfill&\qquad
V_{11}=\{1,5\},\hfill&\b_{11}^{\pm} = \s_4^{\pm 1},\hfill\\
V_6=\{4,5,8\},\hfill&\b_{6}^{\pm}=\s_2^{\pm 1}\s_5^{\pm 1},\hfill&\qquad
V_{12}=\{1,2,3\}.\hfill
\endmatrix
\end{equation*}

An argument as in \cite{CS3} shows that the braid monodromies 
associated to $\W^+$ and $\W^-$ are equivalent, 
but not braid equivalent.  Calculations with these monodromies 
reveal that $\theta^{\loc}_3(\A^\pm) = 16$, 
$\theta_3(\A^\pm) = 21$, and  
$\theta_k(\A^\pm) = \theta^\loc_k(\A^\pm) = 8(k-1)$ for $k\ge 4$. 
Thus the failure of $\overline\Psi_3$
to be surjective is detected only by the third Chen group.  
\end{exm}

\begin{exm}[$9_3$ Configurations]\label{exm:pappus}
The relationship between the (completion of the) Alexander invariant
and the combinatorics of an arrangement appears to be quite delicate.
As an illustration, consider the arrangements $\mathcal P_1$ and 
$\mathcal P_2$ defined by
\begin{align*}
Q(\mathcal P_1)&=xyz(x-y)(y-z)(x-y-z)(2x+y+z)(2x+y-z)(2x-5y+z),\\
Q(\mathcal P_2)&=xyz(x+y)(y+z)(x+3z)(x+2y+z)(x+2y+3z)(2x+3y+3z).
\end{align*}
The arrangement $\mathcal P_1$ is a realization of the Pappus 
configuration $(9_3)_1$, while $\mathcal P_2$ is a realization 
of the configuration $(9_3)_2$.  Note that neither of these
arrangements has a subarrangement lattice-isomorphic to the braid
arrangement.

The combinatorial distinction between these arrangements
(resp.,~their underlying matroids) is detected by the maps
$\overline\Psi_3(\mathcal P_k):\Z^{84}\to\Z^{63}$.  The map 
$\overline\Psi_3(\mathcal P_2)$ is surjective, and consequently
the module $\hB(\mathcal P_2)$ decomposes as a direct sum. Thus, 
$\theta_k(\mathcal P_2)=9 (k-1)$ for $k\ge 2$.  
However, the map $\overline\Psi_3(\mathcal P_1)$ is not surjective, 
and $\hB(\mathcal P_1)$ does not decompose.  A calculation 
shows that 
$\theta_2(\mathcal P_1)=9$, and 
$\theta_k(\mathcal P_1) = 10 (k-1)$ for $k\ge 3$.  

It was conjectured in \cite{CS1} that, for $k$ sufficiently large, 
one has 
$\theta_k(\A)=\theta_k^\loc(\A) + (k-1)\beta(\A)$, 
where $\beta(\A)$ is the number of subarrangements of $\A$ that are 
lattice-isomorphic to $\A_4$.  The arrangement $\mathcal P_1$ has
$\beta=0$ and $\theta_k^\loc = 9(k-1)$, 
and hence provides a counterexample to that conjecture. 
\end{exm}

\bibliographystyle{amsalpha}

\end{document}